\documentclass[11pt,hyp,]{nyjm}
\usepackage{hyperref}
\hypersetup{nesting=true,debug=true,naturalnames=true}
\usepackage{graphicx,amssymb,upref}

\let\<\langle
\let\>\rangle

\let\uml\"

\title{Residual finiteness of $A_{2,3,2n}$ triangle Artin groups}

\author{Greyson Meyer}
\address{Department of Mathematics, University of California Santa Cruz, 1156 High St, Santa Cruz CA 95060 USA.} 
\email{greyson.meyer@gmail.com}  


\usepackage{graphicx}
\usepackage{amssymb}
\usepackage{amsthm}
\usepackage{xcolor}
\usepackage{multirow}
\usepackage{mathtools}
\usepackage{tikz-cd}
\tikzset{
math to/.tip={Glyph[glyph math command=rightarrow]},
loop/.tip={Glyph[glyph math command=looparrowleft, swap]},
loop'/.tip={Glyph[glyph math command=looparrowleft]},
 weird/.tip={Glyph[glyph math command=Rrightarrow, glyph length=1.5ex]},
  pi/.tip={Glyph[glyph math command=pi, glyph length=1.5ex, glyph axis=0pt]},
}
\renewenvironment{proof}{{\bfseries Proof.}}{\qedsymbol}
\newtheorem{theorem}{Theorem}[section]

\newtheorem{lemma}[theorem]{Lemma}
\theoremstyle{definition}
\newtheorem{definition}[theorem]{Definition}

\newtheorem*{remark}{Remark}
\renewcommand\qedsymbol{\raisebox{-4pt}{\includegraphics[height=10pt]{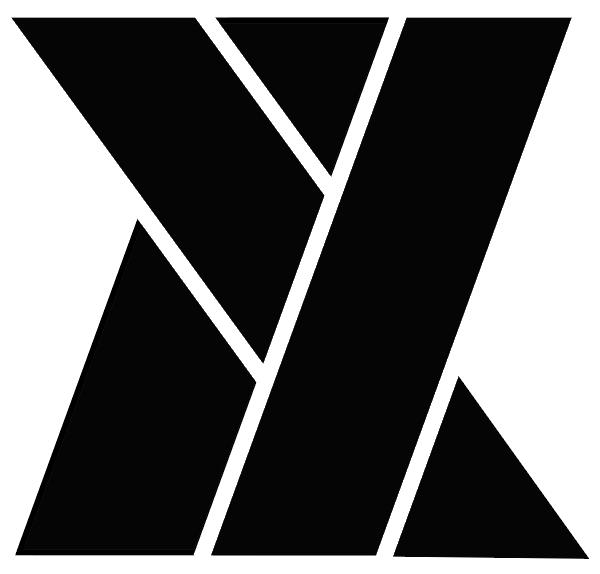}}}

\begin{document} 
 
\begin{abstract}
    We prove that triangle Artin groups of the type $A_{2,3,2n}$ are residually finite for all $n>3$. This requires splitting these triangle Artin groups as graphs of groups and then proving that each of these graphs of groups has finite stature with respect to its vertex groups.
\end{abstract}
\maketitle
\tableofcontents

\section{Introduction}
A \emph{triangle Artin group} is a group that admits the presentation \[A_{M,N,P}=\langle a,b,c|(a,b)_M=(b,a)_M,(b,c)_N=(c,b)_N,(c,a)_P=(a,c)_P\rangle\] where, for example, $(a,b)_M$ is the word with alternating letters in $\{ a,b\}$ with length $M$, beginning with the letter $a$. In \cite{jankiewicz2021splittings}, Jankiewicz proves that every Artin group $A_{M,N,P}$ with $M\leq N\leq P$ and either $M>2$ or $N>3$ splits as a graph of free groups. Jankiewicz then uses this information in \cite{jankiewicz2022residual} and \cite{jankiewicz2023finite} to prove that these same triangle Artin groups are residually finite. A group $G$ is \emph{residually finite} if for every nontrivial $g\in G$, there exists a finite index subgroup $H$ with $g\notin H$. 

The results in \cite{jankiewicz2023finite} exclude two classes of triangle Artin groups, namely $\{A_{2,2,P}|P>1\}$ and $\{A_{2,3,P}|P\geq3\}$. All Artin groups $A_{2,2,P}$ and $\{A_{2,3,P}|3\leq P<6\}$ are spherical, making them linear \cite{cohen2000linearityartingroupsfinite,DIGNE200339} and therefore residually finite \cite{Mal40}. Squier proves in \cite{f4627749-cce2-3922-b300-7d6c84e2c583} that $A_{2,3,6}$ splits as a graph of free groups using purely algebraic means and it can be deduced easily from this result that $A_{2,3,6}$ is residually finite. What remains is the subclass of triangle Artin groups $A_{2,3,P}$ for $P\geq7$. In \cite{wu2023splittings}, Wu \& Ye prove that each $A_{2,3,2p+1}$ with $p\geq3$ cannot split as a graph of free groups. In this paper we will use the splitting found in \cite{wu2023splittings} and expand the methods from \cite{jankiewicz2023finite} to prove the following theorem.
\begin{theorem}
    Every triangle Artin group $A_{2,3,2n}$ with $n>3$ is residually finite.
\end{theorem}

\subsection{Proof strategy}

In \cite{soton50613}, Hanham discovered the following presentation for $A_{2,3,2n}$ with $n>3$, which Wu \& Ye later split in \cite{wu2023splittings} as an amalgamated product of free groups. We condense these two results into Theorem 1.2. For completeness, Theorem 1.2 is proven in Section 2.

\begin{theorem}[Hanham \cite{soton50613}, Wu-Ye \cite{wu2023splittings}]
     When $n>3$,\[
     A_{2,3,2n}
    \cong \langle b,c,x,y,d,\delta|d=xc,d=bx,y=bc,yb=cy,\delta b=c\delta,\delta=d x^{n-2}d\rangle\]
    $\cong F_3*_{F_7}F_4$.
\end{theorem}
The amalgamated product in Theorem 1.2 can be written as $\pi_1(\Gamma)$ where $\Gamma$ is the graph of groups below.
\begin{center}
    \includegraphics[scale=0.15]{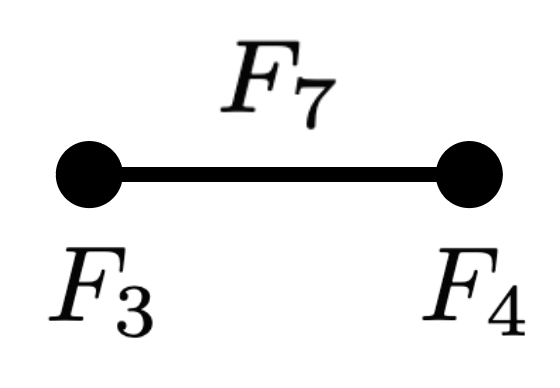}
\end{center}

We will prove that this graph of groups has finite stature with respect to its vertex groups in the sense of \cite{huang2019stature}. We define the concept of finite stature precisely in Section 3. We then will use \cite[Theorem 1.3]{huang2019stature} to deduce residual finiteness. Since this graph of groups is a graph of free groups, the vertex and edge groups implicitly satisfy the first two criteria in \cite[Theorem 1.3]{huang2019stature}, thereby reducing \cite[Theorem 1.3]{huang2019stature} to the more approachable theorem below.

\begin{theorem}[Huang-Wise \cite{huang2019stature}]
    Let $G$ be the fundamental group of a graph of finite rank free groups with finite underlying graph. If $G$ has finite stature with respect to its vertex groups, then $G$ is residually finite.
\end{theorem}

In Section 2 we split each Artin group in $\{A_{2,3,2n}|n>3\}$ as a graph of free groups in a similar manner as in \cite{wu2023splittings}. In Section 3 we develop a procedure involving iterated fiber products used in Sections 4 \& 5 to show that Theorem 1.3 holds in the cases $n>4$ and $n=4$ respectively, proving Theorem 1.1. The proof in Section 5 utilizes a Python program designed by the author to display and categorize the connected components that arise from a multitude of fiber product computations. A link to the GitHub repository containing the code, as well as documentation about how to use the program, can be found in Section 6.

\section{Splitting $A_{2,3,2n}$}
We begin by fixing the presentation of $A_{2,3,2n}$ from Theorem 1.2, which will be used to split $A_{2,3,2n}$ as a graph of free groups. An isomorphism between the standard presentation and the presentation in Theorem 1.2 can be defined by mapping, in one direction:
\[\phi(b)=b, \phi(c)=c,\phi(x)=b^{-1}acb\]
\[\phi(y)=bc, \phi(d)=acb,\phi(\delta)=(ac)^nbc\]
and in the other direction:
\[\psi(a)=db^{-1}c^{-1}, \psi(b)=b,\psi(c)=c\]

The presentation complex associated to the presentation in Theorem 1.2 can be seen in Figure 1. The relations are represented by relator polygons whose sides are identified with the loops in the wedge of circles at the center of the figure. This presentation complex, denoted by $X$, has been constructed so that $\pi_1(X)= A_{2,3,2n}$, endowed with the presentation from Theorem 1.2. In order to realize $A_{2,3,2n}$ as an amalgamated product, we apply Seifert-Van Kampen's Theorem to $X$ as shown in Figure 2.
\begin{figure}[h]
    \centering
    \includegraphics[scale=0.25]{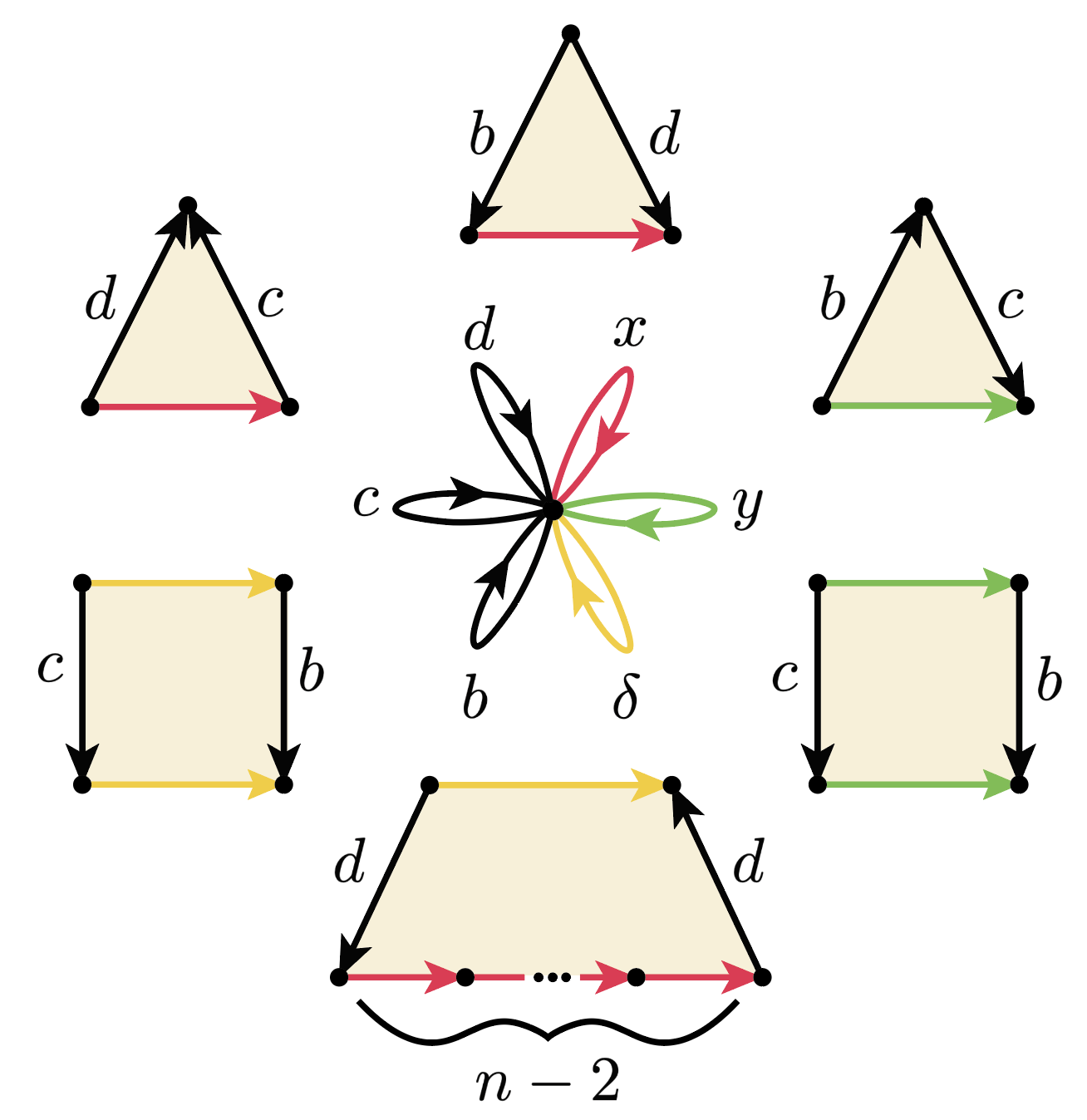}
    \caption{Presentation complex $X$}
    \label{fig:enter-label}
\end{figure}
\begin{figure}[h]
    \centering
    \includegraphics[scale=0.25]{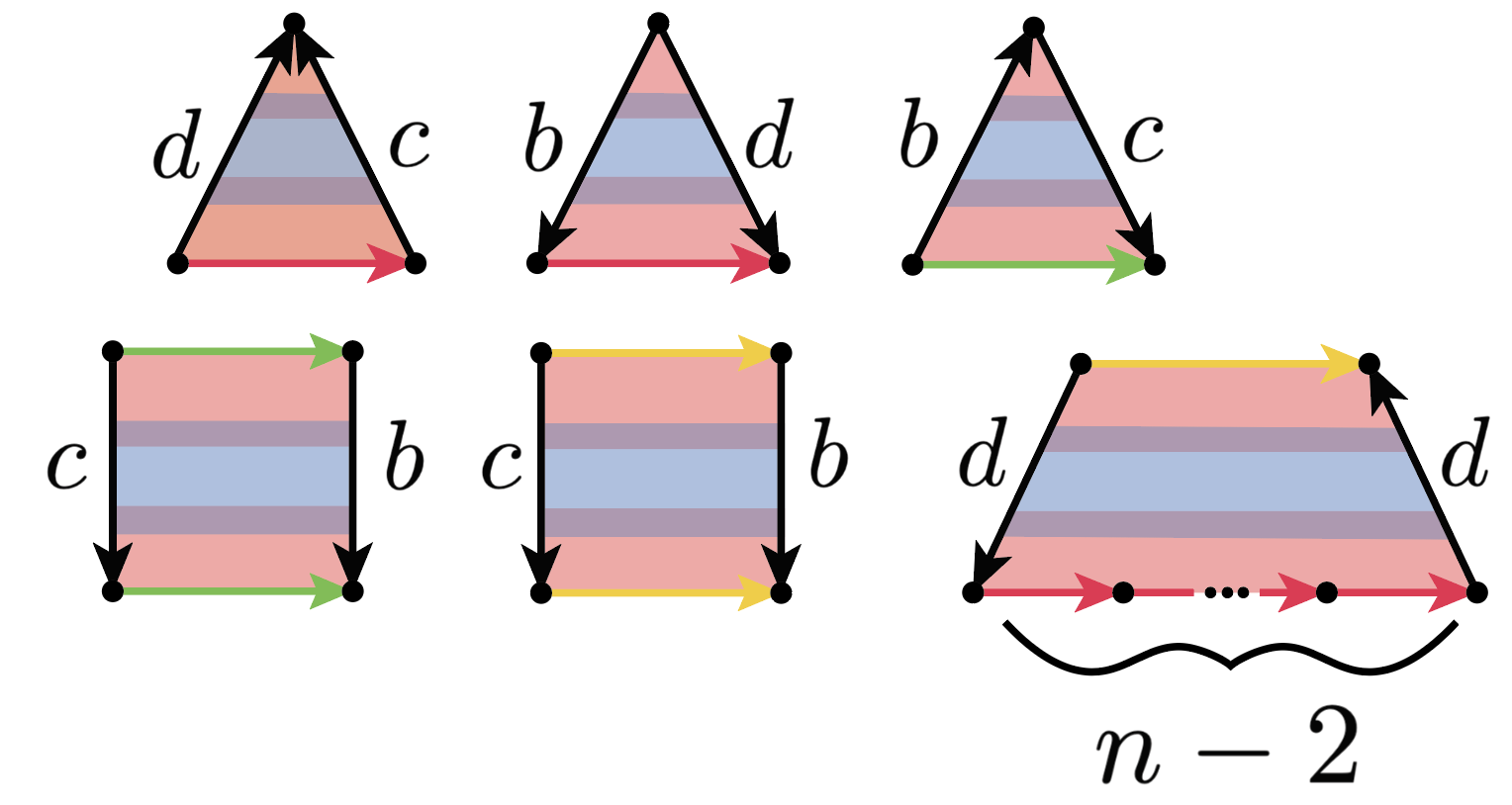}
    \caption{Seifert Van Kampen's Theorem applied to $X$}
    \label{fig:enter-label}
\end{figure}

The red band $N_U$, the blue band $N_V$ and their nonempty intersection, the purple band $N_W$, are path connected open subsets of $X$. There is a natural deformation retraction that we then perform on $N_U$ and $N_V$, respectively, to obtain the graphs $U$ and $V$, as shown in Figures 3 \& 4.
\begin{figure}[h!]
    \centering
    \includegraphics[scale=0.31]{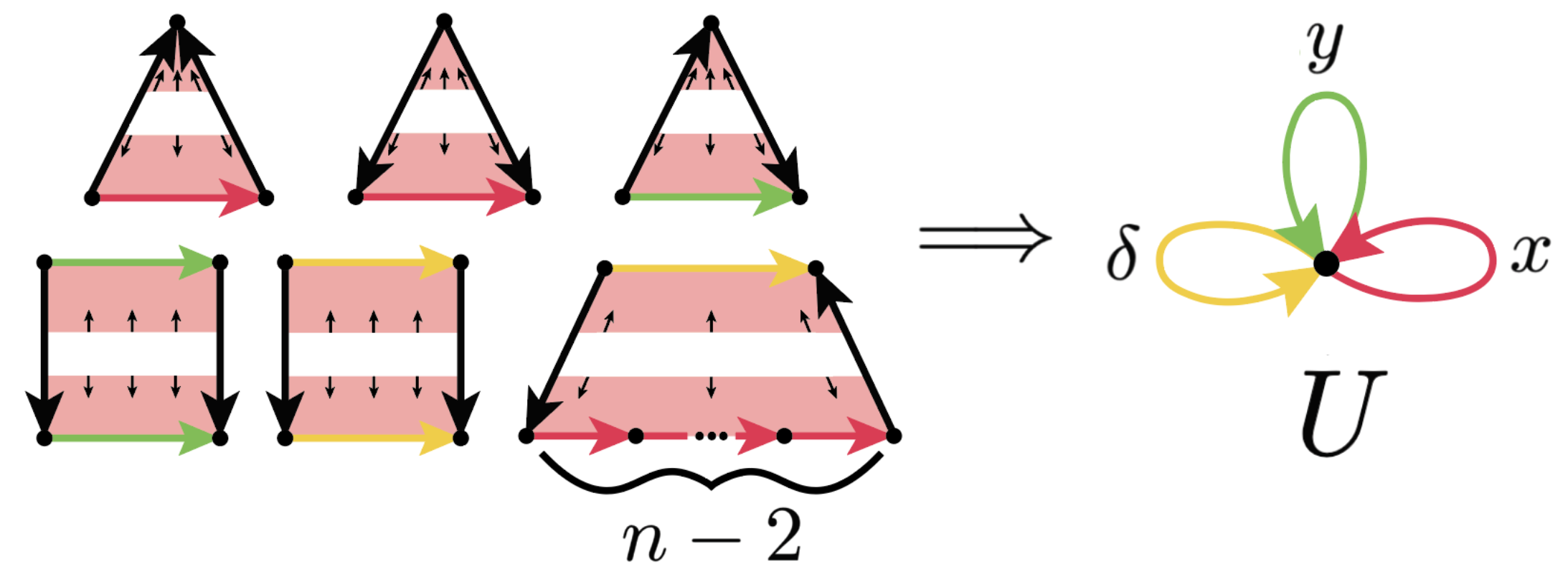}
    \caption{$N_U$ deformation retracts to the wedge of 3 circles, $U$.}
    \label{fig:enter-label}
\end{figure}
\begin{figure}[h]
    \centering
    \includegraphics[scale=0.3]{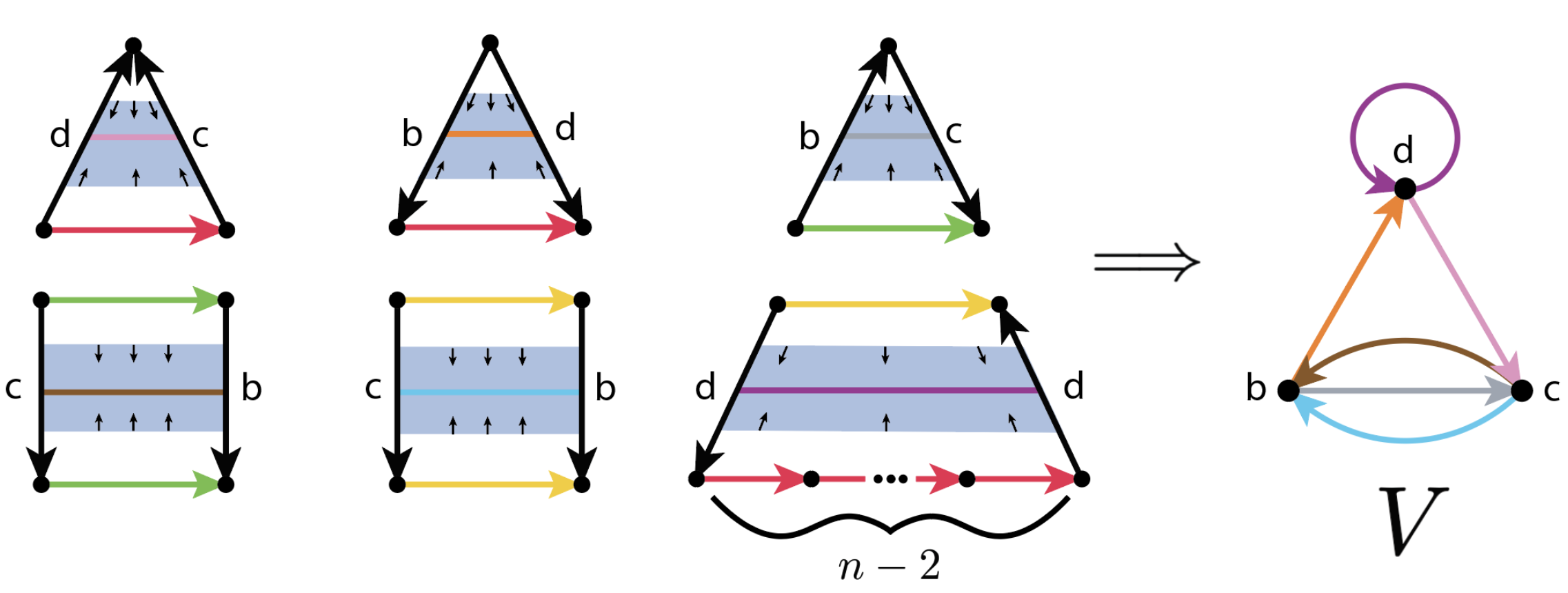}
    \caption{$N_V$ deformation retracts to the labeled graph $V$ with $\pi_1(V)\cong F_4$.}
    \label{fig:enter-label}
\end{figure}

Next, we deformation retract $N_W$ to the labeled graph $W$ represented by the thin purple horizontal lines that appear twice in each relator polygon in Figures 5 \& 6. Note that there is a natural map $W\to V$ induced by inclusions and the deformation retractions. We
will soon see that this map is a combinatorial immersion. 

\begin{definition}[\cite{Stallings1983}]
A combinatorial map $Y\to Z$ between graphs $Z$ and $Y$ is a function that maps vertices to vertices and edges to edges. A combinatorial immersion $\phi:Y\looparrowright Z$ between these graphs is a locally injective combinatorial map.
\end{definition}

Every combinatorial immersion $\phi:Y\looparrowright Z$ induces an injective homomorphism $\pi_1(Y,y)\hookrightarrow\pi_1(Z,z)$ \cite{Stallings1983}. Furthermore, if $Y$ is not already a cover of $Z$, the existence of a combinatorial immersion $\phi:Y\looparrowright Z$ guarantees that $Y$ can be completed to a cover of $Z$ by attaching trees to the vertices in $Y$ that are inhibiting $Y$ from being a cover  \cite{Stallings1983}.

\begin{lemma}
    The map $W\looparrowright V$ induced by the deformation retraction $N_V\to V$ is a combinatorial immersion.
\end{lemma}
\begin{proof}
The map $W\looparrowright V$ factors through the following composition of maps:
\[W\hookrightarrow N_W\hookrightarrow N_V\to V\]
where the rightmost map is the deformation retraction of $N_V$. The vertex labels in $W$ come directly from the relator polygons. The sides of the relator polygons are labeled $b$, $c$ or $d$. Every vertex in $W$ lies on one such edge. The superscript is then included to keep track of whether the vertex is closer to the head (in which case the superscript is $+$) or the tail (in which case the superscript is a $-$) of the oriented edge. The combinatorial immersion is represented by the coloring of the edges in $W$.
\end{proof}

\begin{figure}[h]
    \centering
    \includegraphics[scale=0.38]{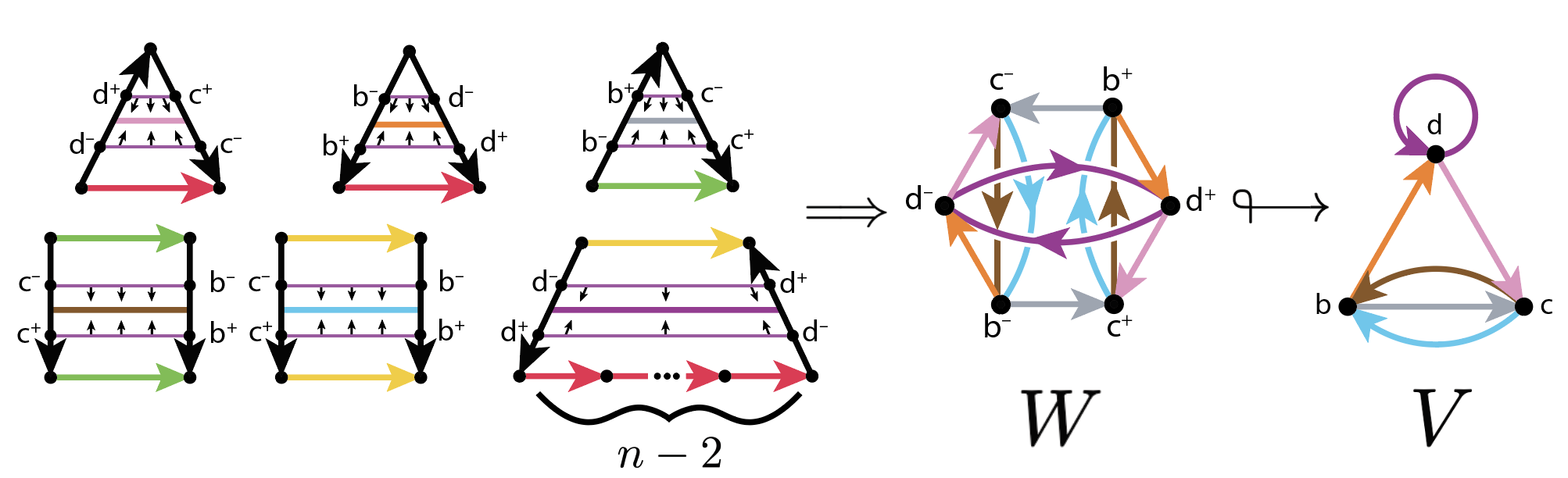}
    \caption{The labeled graph $W$ is a double cover of $V$.}
    \label{fig:enter-label}
\end{figure}

Notice further that $W$ is a double cover of $V$. The map $W\to U$ induced by the deformation retraction $N_U\to U$ is not as simple. Just like in Lemma 2.2, the map $W\to U$ factors through the composition of maps
\[W\hookrightarrow N_W\hookrightarrow N_U\to U\]
where the rightmost map is the deformation retraction. What is keeping the map $W\to U$ from being a combinatorial immersion is the triangular relator polygons. The deformation retraction $N_U\to U$ collapses the topmost internal edges of the triangular relator polygons to vertices. To take these collapsing edges into account, we denote the graph in the bottom right of Figure 6 as $\sigma(W)$ and decompose the map $W\to U$ as the composition $W\xrightarrow{\sigma}\sigma(W)\looparrowright U$ where the second map is a combinatorial immersion. The map $\sigma$ and the combinatorial immersion are shown in Figure 6. Again, the combinatorial immersion is represented by the colors of the edges in $\sigma(W)$.

\begin{lemma}
    The map $\sigma:W\to\sigma(W)$ is a homotopy equivalence.
\end{lemma}
\begin{proof}
    Every edge in $V$ arises from a unique relator polygon. The relator polygon from which each edge originates is encoded in the color of the edge. Intuitively, $\sigma$ sends every edge in $W$ to the nearest horizontal edge of the relator polygon from which that edge originates. To see how $\sigma$ maps every edge in $W$, see Figure 6. Mostly $\sigma$ sends an edge in $W$ to an edge in $\sigma(W)$. The exceptions are the edges $d^+\to c^+$, $b^-\to d^-$ and $b^+\to c^-$ in $W$ that are collapsed to vertices (represented by the dashed black lines in Figure 6), and the edge $d^+\to d^-$ in $W$ that is sent to the path of length $n-2$ of red edges (the curve labeled ``$n-2$" in Figure 6). Since arcs are contractible, all four of these nontrivial steps that comprise the map $\sigma$ are homotopy equivalences, making $\sigma$ itself a homotopy equivalence as well.
\end{proof}

While $W$ is a double cover of $V$, $\sigma(W)$ is not a cover of $U$. The red arc at the bottom of the graphs in Figure 6 labeled ``$n-2$" represents an oriented path consisting of $n-2$ red edges. The vertices in this path are not the endpoints of any green or yellow edges, making these vertices not preimages of the lone vertex in $U$. Graphs of this form can be completed to a cover by attaching trees to the vertices that lack the proper adjacencies. In the case of $\sigma(W)$, this requires attaching an appropriate subtree of the universal cover of $U$ to each vertex in the red path since the other three vertices in $\sigma(W)$ have full adjacency. This act of completing a graph to a cover by attaching trees does not affect the fundamental group of the graph since trees are devoid of loops. We refer to $\sigma(W)$ as the \emph{core} of a cover of $U$, which is the portion of the cover containing loops \cite{Stallings1983}. 

\begin{figure}[h!]
    \centering
    \includegraphics[scale=0.35]{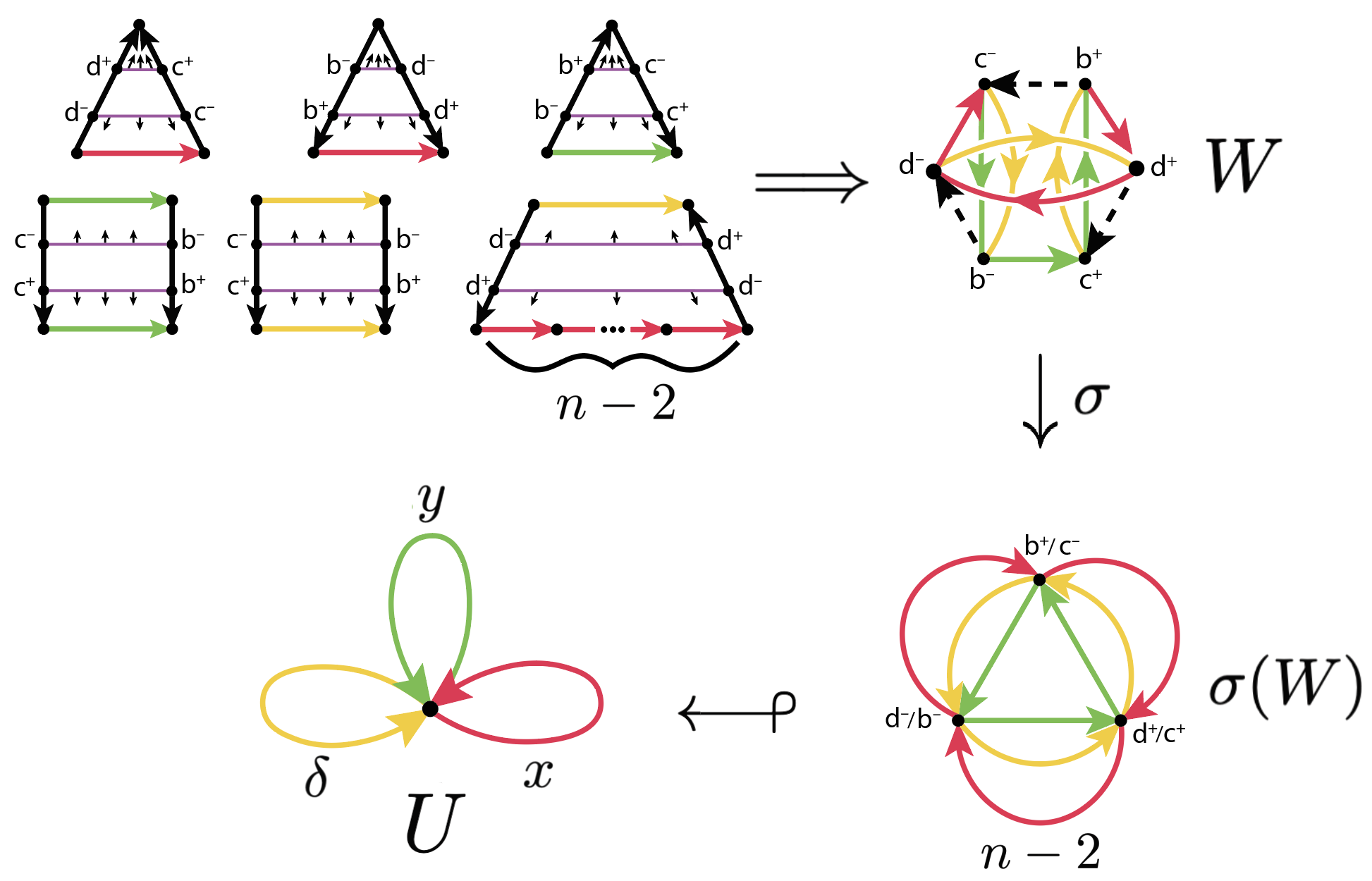}
    \caption{The labeled graph $\sigma(W)$ combinatorially immerses into $U$.}
    \label{fig:enter-label}
\end{figure}

We use the combinatorial immersions $\sigma(W)\looparrowright U$ and $W\looparrowright V$ to induce the injections $\pi_1(\sigma(W))\hookrightarrow\pi_1(U)$ and $\pi_1(W)\hookrightarrow\pi_1(V)$ respectively. Since we would need infinitely many vertices to complete $\sigma(W)$ to a cover of $U$ with the same fundamental group as $U$, and $U$ has only one vertex, this makes $\pi_1(\sigma(W))$ an infinite-index subgroup of $\pi_1(U)$. Since $W$ is a double cover of $V$, this makes $\pi_1(W)$, an index-2 subgroup of $\pi_1(V)$. The deformation retractions, combinatorial immersions and $\sigma$ all fit into the following diagram that commutes up to homotopy:

\begin{center}
\begin{tikzcd}
\sigma(W) \arrow[d, loop-math to] \\
U \arrow[d, hook] & W \arrow[d, hook] \arrow[r, loop-math to] \arrow[ul, "\sigma"]
  & V \arrow[d, hook] \\
N_U & N_W \arrow[l] \arrow[r] & N_V
\end{tikzcd}
\end{center}

This diagram, along with Seifert-Van Kampen's Theorem being applied to the presentation complex for $A_{2,3,2n}$, proves Theorem 2.4 in the same manner as in \cite{wu2023splittings}.

\begin{theorem}\cite{wu2023splittings}
    For $n\geq3$, $A_{2,3,2n}\cong \pi_1(N_U)*_{\pi_1(N_W)}\pi_1(N_V)$ where \begin{itemize}
    \item $\pi_1(N_U)=\pi_1(U)\cong F_3$
    \item $\pi_1(N_V)=\pi_1(V)\cong F_4$
    \item $\pi_1(N_W)=\pi_1(W)=\pi_1(\sigma(W))\cong F_7$
\end{itemize}
\end{theorem}

\section{Finite stature procedure}
In order to prove that $A_{2,3,2n}$ is residually finite, we first prove that $A_{2,3,2n}$ has finite stature with respect to $\{\pi_1(U),\pi_1(V)\}$, the vertex groups of the splitting from Theorem 2.4. The definition of finite stature can be found in Definition 3.6 in Section 3.3. The definition involves intersections of subgroups and conjugacy classes of subgroups. So far we have only discussed the groups in our graph of groups in terms of the fundamental group of labeled graphs. We would like to continue to proceed in this fashion. In order to approach the concept of finite stature from this perspective, we will need to understand what subgroup intersection and conjugation look like in the context of labeled graphs.

\subsection{Fiber products}
Stallings describes in \cite{Stallings1983} how to calculate the intersection of subgroups that can be realized as the fundamental group of respective covers of a common underlying labeled graph. He does so using the following tools.

\begin{definition}
Let $\phi_i:Y_i\to X$ be a combinatorial immersion for $i=1,2$. The \emph{fiber product} of $Y_1$ and $Y_2$ over $X$ is the labeled graph $Y_1\otimes_XY_2$ with vertex set\[
\{(v_1,v_2)\in V(Y_1)\times V(Y_2):\phi_1(v_1)=\phi_2(v_2)\}\] and edge set\[
\{(e_1,e_2)\in E(Y_1)\times E(Y_2):\phi_1(e_1)=\phi_2(e_2)\}\]
There is a natural combinatorial immersion $Y_1\otimes_XY_2\to X$, given by $(y_1,y_2)\mapsto\phi_1(y_1)=\phi_2(y_2)$.
\end{definition}
\begin{lemma}[Stallings \cite{Stallings1983}]
    Let $H_1,H_2\leq \pi_1(X,v)$ where $X$ is a finite graph and $v\in X$ is a vertex. For $i=1,2$, let $(Y_i,\hat{x}_i)\to(X,v)$ be a cover of $X$ where $\pi_1(Y_i,\hat{x}_i)= H_i$. Then $H_1\cap H_2=\pi_1(Y_1\otimes_XY_2,(\hat{x}_1,\hat{x}_2))$.
\end{lemma}

We will be performing fiber products of covers of $W$ and $\sigma(W)$ throughout the rest of the paper. Since $W$ and $\sigma(W)$ are homotopy equivalent, we choose to work with covers of $\sigma(W)$ instead of $W$. We make this choice since $U$ has only one vertex, which guarantees that the immersions $\phi_1$ and $\phi_2$ in Definition 3.1 will agree on every vertex. In terms of Definition 3.1 this means that $U$ will be the graph that we will perform fiber products over (the $X$ in Definition 3.1) and every labeled graph $Y_i$ will combinatorially immerse into $U$. This will result in $V(Y_1)\times V(Y_2)$ being the vertex set of every fiber product for the rest of this paper for all $Y_1,Y_2\looparrowright U$. Figure 7 shows an example of a fiber product calculation that will be used later on.

\begin{figure}[h]
    \centering
    \includegraphics[scale=0.35]{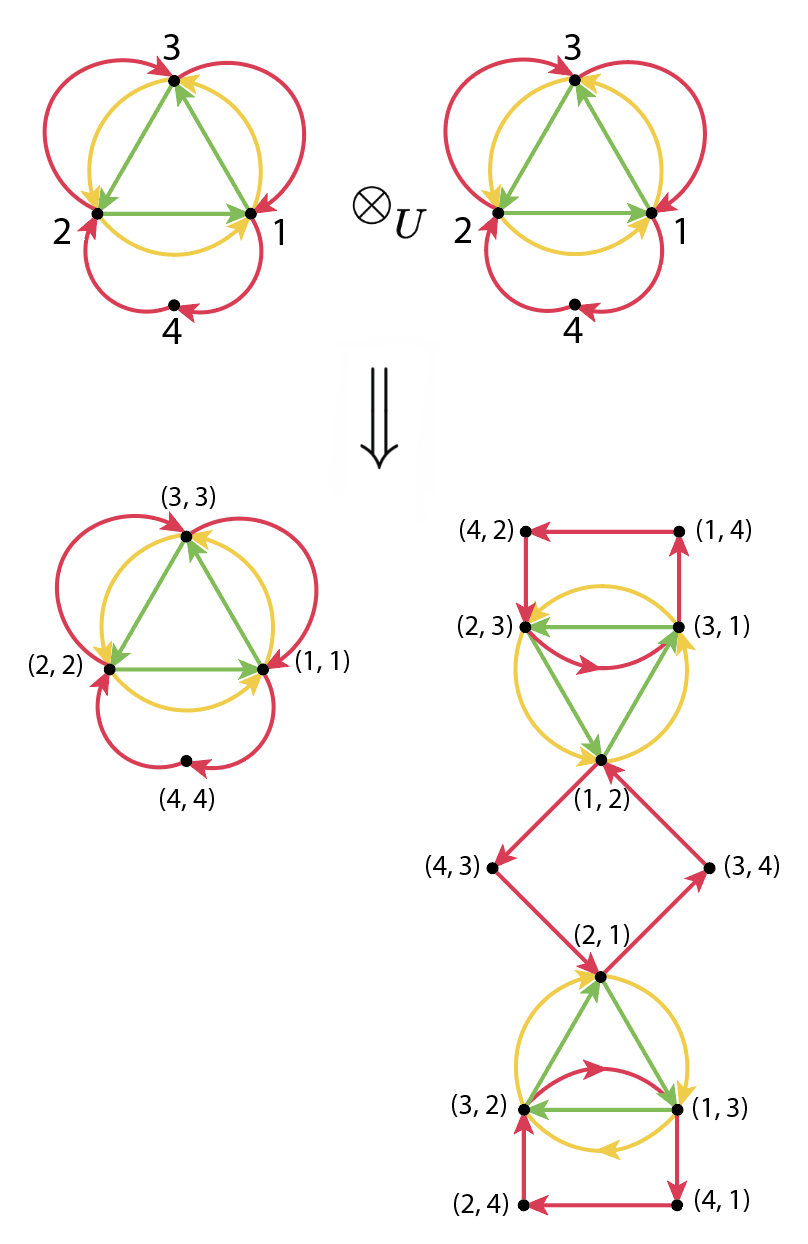}
    \caption{$\sigma(W)\otimes_U\sigma(W)$ for $A_{2,3,8}$}
    \label{fig:enter-label}
\end{figure}

Notice that every element of $V(\sigma(W))\times V(\sigma(W))$ is present in the $\sigma(W)\otimes_{U}\sigma(W)$ calculation in Figure 7, where we have replaced the vertex labels of $\sigma(W)$ with numbers to keep notation clean. We will continue with this labeling convention for the rest of the paper. Figure 7 also demonstrates that fiber products need not be connected, making choice of basepoint very important when calculating the fundamental group. Choosing basepoints $a\in Y_1$ and $b\in Y_2$, makes $(a,b)$ the basepoint in $Y_1\otimes_{U} Y_2$. We will still be primarily interested in the cores of graphs, and thus will continue to omit the trees needed to complete the cores that appear as connected components of fiber products to actual covers of $U$. The following lemma will come in handy during future calculations.

\begin{lemma}
    Let $Y_1\subset X_1$ and $Y_2\subset X_2$ be subgraphs. Then $Y_1\otimes_{U} Y_2\subset X_1\otimes_{U} X_2$.
\end{lemma}
This lemma is immediate by the definition of fiber products. Also, from an algebraic point of view, the intersection of subgroups will always be a subgroup of the intersection of the supergroups.

\subsection{Conjugation of subgroups as operations on labeled graphs}

We now must discuss how to understand conjugation of subgroups when our subgroups are the fundamental group of labeled graphs. To start, we will fix $v\in\pi_1(V)$ as being a fixed coset representative coming from $\pi_1(V)/\pi_1(W) = \{[1],[v]\}$. This coset representative allows us to observe the following property of $A_{2,3,2n}$ for $n\geq3$. 

\begin{lemma}
    Every element $g\in A_{2,3,2n}$ with $n\geq3$ can be written as \[g=g_1vg_2v...g_{m-1}vg_m\] for some $m\geq1$ where $g_1,...,g_m\in\pi_1(U)$.
\end{lemma}
\begin{proof}
The splitting in Theorem 2.4 gives us \[A_{2,3,2n}\cong\pi_1(U)*_{\pi_1(W)=\pi_1(\sigma(W))}\pi_1(V)=(\pi_1(U)*\pi_1(V))/(\pi_1(W)=\pi_1(\sigma(W)))\]
Therefore every element $g$ can be written $g=\bar{h}_1\bar{k}_1\bar{h}_2\bar{k}_2...\bar{h}_{\mathfrak{m}-1}\bar{k}_{\mathfrak{m}-1}\bar{h}_\mathfrak{m}$ for some finite $\mathfrak{m}$ and where each $\bar{h}_i\in\pi_1(U)$ and $\bar{k}_i\in\pi_1(V)$. Since $\pi_1(W)$ is an index-2 subgroup of $\pi_1(V)$, each $\bar{k}_i=l_i$ or $\bar{k}_i=l_iv$ for some $l_i\in\pi_1(\sigma(W))=\pi_1(W)$. In the first case $\bar{h}_i\bar{k}_i\bar{h}_{i+1}\in\pi_1(U)$, allowing us to replace $\bar{h}_i\bar{k}_i\bar{h}_{i+1}$ with some $h_i\in\pi_1(U)$ to rewrite $g=h_1k_1...k_{m-1}h_m$ where $m\leq \mathfrak{m}$ and every $k_i$ is of the form $l_iv$. Therefore $g=h_1k_1...h_{m-1}k_{m-1}h_m=h_1l_1vh_2l_2v...h_{m-1}l_{m-1}vh_m$. Setting $h_il_i=g_i$ for $1<i<m$ and $h_m=g_m$ results in the desired expansion of $g$.
\end{proof}

Let $p_u$ and $p_\sigma$ be fixed basepoints of $U$ and $\sigma(W)$ respectively. Consider a cover $(Y,p)\to(\sigma(W),p_\sigma)$ and let $\pi_1(Y,p)=H\leq\pi_1(\sigma(W),p_\sigma)$ be the corresponding subgroup. By Lemma 3.4 it suffices to understand how conjugation by elements of $\pi_1(U)$, and conjugation by the coset representative $v\in\pi_1(V)$, affect $H$.

Since $\sigma(W)$ is a cover of $U$, this makes $(Y,p)\to(U,p_u)$ a cover as well. So the conjugation of $H$ by an element of $u\in U$ in terms of the graph $Y$ is represented by simply shifting $p$ along the path in $Y$ described by $u$ (see \cite{MR1867354} for more details on the connection between basepoint shifting and conjugation).

But $Y$ is not a cover of $W$, nor of $V$, so we cannot represent the $v$-conjugation of $H$ by simply shifting $p$ along a path in $Y$. Since $W$ is a double cover of $V$, $v$ is a path between the two preimages in $W$ of the basepoint in $V$ (see Figure 5 for the vertex labellings of $V$ and $W$). This takes the basepoint of $W$ from one preimage of $V$ to the other, resulting in every vertex label having its superscript swapped, meaning
\begin{center}
$d^+\leftrightarrow d^-,b^+\leftrightarrow b^-$ and $c^+\leftrightarrow c^-$
\end{center}

This involution extends naturally to all labeled graphs that combinatorially immerse into $W$.

\begin{definition}
    Define the involution $\bar{\beta}:\{Y'\looparrowright W\}\to\{Y'\looparrowright W\}$ on the vertices of a labeled graph $Y'$ by $d^+\leftrightarrow d^-$, $b^+\leftrightarrow b^-$ and $c^+\leftrightarrow c^-$.
\end{definition} 

Consider the homotopy equivalence $\sigma:W\to \sigma(W)$ introduced in Section 2. Since $\sigma$ is a homotopy equivalence, it extends naturally to all covers of $W$ as well. Figure 8 shows how $\sigma$ behaves locally on edges. Given some labeled $Y\looparrowright \sigma(W)$, if we want to find a labeled graph $Z\looparrowright \sigma(W)$ such that $\pi_1(Z, *_v)=v^{-1}\pi_1(Y, *)v$, we start by considering the labeled graph $Y'\looparrowright W$ such that $\sigma(Y')=Y$. The $v$-translated copy of $Y'$ is the labeled graph $\bar{\beta}(Y')$. This new labeled graph $\bar{\beta}(Y')$ arises as a copy of $Y'$ with $+\leftrightarrow -$ in all of the vertex label superscripts. 

Putting these pieces together, we define $\beta:\{Y\looparrowright\sigma(W)\}\to\{Y\looparrowright \sigma(W)\}$ by $\beta(Y)=\sigma\circ\bar{\beta}(Y')$ where $\sigma(Y')=Y$. Note that $\beta^2=1$ by construction. Pictorially, we can think of $\beta(\sigma(W))$ as swapping the tops and bottoms of the relator polygons in the presentation complex. Figure 8 shows the details of how $\beta$ affects the vertices and edges of the core of every cover of $\sigma(W)$. Note that since $\beta$ is an involution, Figure 8 is meant to also be read from right to left. In the figure, the dashed arrows are meant to represent edges that are being mapped to a vertex by $\sigma$. An example of a $\beta$ calculation can be found in Figure 12.

\begin{figure}[h]
    \centering
    \includegraphics[scale=0.45]{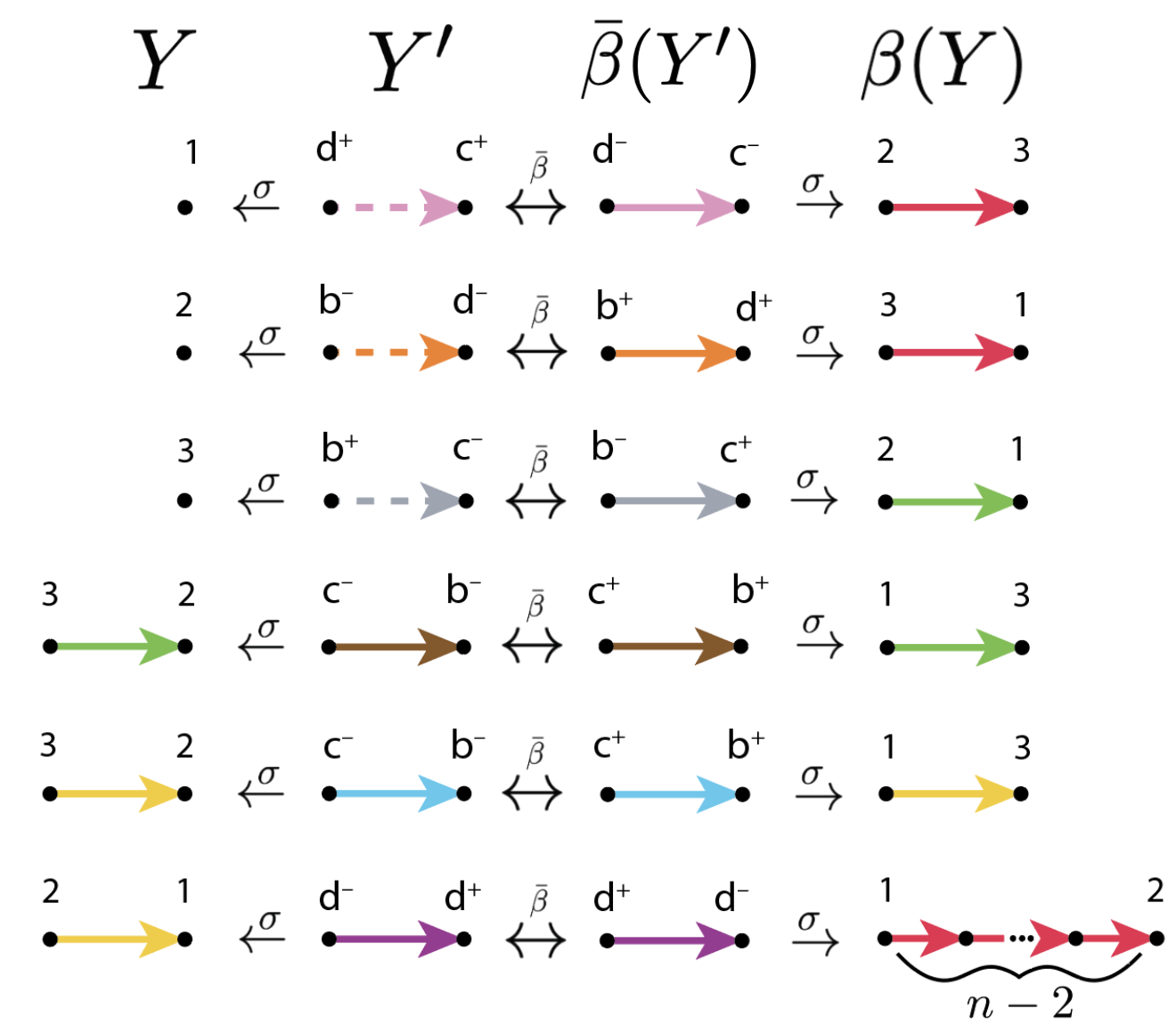}
    \caption{$\beta(Y)=\sigma\circ\bar{\beta}(Y')$ applied to every applicable vertex \& edge.}
    \label{fig:enter-label}
\end{figure}

To recap, for every labeled graph $Y\looparrowright \sigma(W)$, the labeled graph corresponding to $g^{-1}\pi_1(Y,*)g$ for some $g\in A_{2,3,2n}$ is either:
\begin{enumerate}
    \item $Y$ with a new basepoint chosen, when $g\in\pi_1(U)$ or
    \item a (potentially) new labeled graph $\beta(Y)$ with a (potentially) new basepoint, when $g\notin\pi_1(U)$.
    \end{enumerate}
\subsection{Defining the set $S$}

We finally have enough tools to bring finite stature into the mix.

\begin{definition}[{\cite{huang2019stature}}]
    Let $G$ be a group and $\{H_\lambda\}_{\lambda\in\Lambda}$ be a collection of subgroups of $G$. Then $G$ has \emph{finite stature} with respect to $\{H_\lambda\}_{\lambda\in\Lambda}$ if for each $\mu\in\Lambda$, there are finitely many $H_\mu$-conjugacy classes of infinite subgroups of the form $H_\mu\cap D$ where $D$ is an intersection of $G$-conjugates of elements in $\{H_\lambda\}_{\lambda\in\Lambda}$.
\end{definition} 

\begin{remark}
Connected components coming from fiber products have tuples for vertex labels. The maps $\sigma$ and $\beta$ are defined with respect to vertices that have integer labels. We will need to calculate the image under $\beta$ of connected components coming from fiber products, so we choose to project each vertex-tuple to its second component. Because we are making this choice, we must perform $H\otimes_{U}K$ and $K\otimes_{U}H$ for all pairs of finite connected labeled graphs $H,K\looparrowright U$ going forward, despite the fact that $(H\otimes_{U}K)\cong(K\otimes_{U}H)$ is an isomorphism of labeled graphs realized by swapping the tuple-entries at each vertex.
\end{remark}

Lemma 3.9 converts the search for finite stature to the search for a set of labeled graphs $S$ with the following important properties.

\begin{definition}
    Let $S$ be the set of finite connected labeled graphs that combinatorially immerse into $\sigma(W)$, defined by the following recursive definition\begin{itemize}
        \item $\sigma(W)\in S$,
        \item if $X,Y\in S$ and $Z$ is a connected component of $X\otimes_UY$, then $Z\in S$,
        \item if $X\in S$, then $\beta(X)\in S$.
    \end{itemize}
    These rules determine $S$, but $S$ is a priori an infinite set.
\end{definition}

\begin{lemma}
    For all $g\in A_{2,3,2n}$, $\pi_1(U)\cap\pi_1(\sigma(W))^g$ is contained in a $\pi_1(U)$-conjugacy class of subgroups of $\pi_1(U)$ represented by a labeled graph in $S$. Likewise, $\pi_1(V)\cap\pi_1(\sigma(W))^g$ is also contained in a $\pi_1(U)$-conjugacy class of subgroups of $\pi_1(U)$ represented by a labeled graph in $S$.
\end{lemma}
\begin{proof}
    Lemma 3.4 allows us to decompose every $g\in A_{2,3,2n}$ as \[g=g_1vg_2v...g_{m-1}vg_m\]
    where $m\geq1$,  $g_i\in\pi_1(U)$ and $v\in\pi_1(V)$ is the fixed nontrivial coset representative coming from the quotient $\pi_1(V)/\pi_1(W)\cong\mathbb{Z}/2\mathbb{Z}$. Define $\ell(g)$ to be the number of nontrivial elements in the decomposition of $g$ as above. We will proceed via induction on $\ell(g)$.

    Base case: Consider an element $g\in A_{2,3,2n}$ where $\ell(g)=1$. Then $g=v$ or $g$ is some nontrivial element in $\pi_1(U)$.
    
    First assume $g=v$. Since $\pi_1(W)^v=\pi_1(W)=\pi_1(\sigma(W))$, we have $\pi_1(U)\cap\pi_1(\sigma(W))^g=\pi_1(V)\cap\pi_1(\sigma(W))^g=\pi_1(\sigma(W))$.

    Now assume $g\in\pi_1(U)$. Then $\pi_1(\sigma(W))^g<\pi_1(U)$, making $\pi_1(U)\cap\pi_1(\sigma(W))^g=\pi_1(\sigma(W))^g$. The $\pi_1(U)$-conjugacy class of this subgroup is represented by $\sigma(W)$ with its basepoint shifted according to $g$. The labeled graph $\sigma(W)\in S$ by construction.

    In the Bass-Serre tree $T$ associated with the splitting of $A_{2,3,2n}$ arising from Theorem 2.4, we see that the vertex stabilized by $\pi_1(U)$ is between the edge stabilized by $\pi_1(\sigma(W))^g$ and the vertex stabilized by $\pi_1(V)$. Therefore $\pi_1(U)$ is also stabilized by $\pi_1(V)\cap\pi_1(\sigma(W))^g$. So we can rewrite $\pi_1(V)\cap\pi_1(\sigma(W))^g=\pi_1(V)\cap\pi_1(U)\cap\pi_1(\sigma(W))^g=\pi_1(\sigma(W))\cap\pi_1(\sigma(W))^g$. Therefore, by Lemma 3.2, the labeled graph associated with the $\pi_1(U)$-conjugacy class of this subgroup is a connected component of $\sigma(W)\otimes_U\sigma(W)$, all of which are contained in $S$ since $S$ is closed under fiber product.

    Induction hypothesis: Assume that the lemma holds for all $h\in A_{2,3,2n}$ with $\ell(h)<k$ for some $k>1$.
    
    Induction step: Consider an element $g\in A_{2,3,2n}$ with $\ell(g)=k$.
    
    Case 1:
Let $g=g'v$ where the decomposition of $g'$ begins with $v$ and $\ell(g')=k-1$. Consider the Bass-Serre tree $T$ induced by the splitting of $A_{2,3,2n}$. The subgroup $\pi_1(U)\cap\pi_1(\sigma(W))^g$ is the stabilizer of the path $p$ in $T$ between the vertex in $T$ stabilized by $\pi_1(U)$ and the edge in $T$ stabilized by $\pi_1(\sigma(W))^g$. Since $g$ begins with $v$, this makes the vertex stabilized by $\pi_1(U)^v$ on the path $p$. Therefore $\pi_1(U)\cap\pi_1(\sigma(W))^g=\pi_1(U)\cap(\pi_1(U)\cap\pi_1(\sigma(W))^{g'})^v$. By the induction hypothesis, $\pi_1(U)\cap\pi_1(\sigma(W))^{g'}$ is contained in a $\pi_1(U)$-conjugacy class represented by a labeled graph $H^1_U\in S$. The $v$-conjugate of a subgroup contained in the $\pi_1(U)$-conjugacy class represented by $H^1_U$ is contained in the $\pi_1(U)$-conjugacy class represented by $\beta(H^1_U)$ by the definition of $\beta$. The graph $\beta(H^1_U)\in S$ since $S$ is closed under $\beta$.

Similarly, $\pi_1(V)\cap\pi_1(\sigma(W))^g=(\pi_1(V)\cap\pi_1(\sigma(W))^{g'})^v$ since $v\in\pi_1(V)$. By the induction hypothesis, $\pi_1(V)\cap\pi_1(\sigma(W))^{g'}$ is contained in a $\pi_1(U)$-conjugacy class represented by a labeled graph $H^1_V\in S$. Therefore $\pi_1(V)\cap\pi_1(\sigma(W))^g$ is represented by $\beta(H^1_V)\in S$.

Case 2: Now assume that $g=g'v$ and the decomposition of $g'$ begins with some nontrivial $g_1\in\pi_1(U)$. Since $k>1$, this forces $g$ to begin with $g_1v$. The vertex in $T$ stabilized by $\pi_1(U)^{g_1v}=\pi_1(U)^v$ is on the path in $T$ between the vertex in $T$ stabilized by $\pi_1(U)$ and the edge in $T$ stabilized by $\pi_1(\sigma(W))^g$. Therefore $\pi_1(U)\cap\pi_1(\sigma(W))^g=\pi_1(U)\cap\pi_1(U)^{g_1v}\cap\pi_1(\sigma(W))^g=\pi_1(U)\cap(\pi_1(U)^{g_1}\cap\pi_1(\sigma(W))^{g'})^v=\pi_1(U)\cap(\pi_1(U)\cap\pi_1(\sigma(W))^{g'})^v$. By the induction hypothesis, $\pi_1(U)\cap\pi_1(\sigma(W))^{g'}$ is contained in the $\pi_1(U)$-conjugacy class represented by a labeled graph $H^2_U\in S$. So $\pi_1(U)\cap\pi_1(\sigma(W))^g=\pi_1(U)\cap(\pi_1(U)\cap\pi_1(\sigma(W))^{g'})^v$ is represented by $\beta(H^2_U)\in S$.

Similarly, $\pi_1(V)\cap\pi_1(\sigma(W))^g=(\pi_1(V)\cap\pi_1(\sigma(W))^{g'})^v$. By the induction hypothesis we know that $\pi_1(V)\cap\pi_1(\sigma(W))^{g'}$ is contained in the $\pi_1(U)$-conjugacy class represented by a labeled graph $H^2_V\in S$. Therefore $\pi_1(V)\cap\pi_1(\sigma(W))^g$ is represented by $\beta(H^2_V)\in S$.

Case 3: Assume $g=g'g_m$, $g_m\neq1$, the decomposition of $g'$ begins with $v$ and $\ell(g')=k-1$. Then $\pi_1(U)\cap\pi_1(\sigma(W))^g=(\pi_1(U)\cap\pi_1(\sigma(W))^{g'})^{g_m}$ since $g_m\in\pi_1(U)$. By the induction hypothesis we know that $\pi_1(U)\cap\pi_1(\sigma(W))^{g'}$ is contained in the $\pi_1(U)$-conjugacy class represented by a labeled graph $H^3_U\in S$. Therefore $\pi_1(U)\cap\pi_1(\sigma(W))^g$ is also represented by $H^3_U$.

Define $g''$ to be such that $g=vg''$ with $\ell(g'')=k-1$ and $g''$ ends with $g_m$. Then $\pi_1(V)\cap\pi_1(\sigma(W))^g=\pi_1(V)\cap((\pi_1(V)\cap\pi_1(\sigma(W)))^v)^{g''}=\pi_1(V)\cap(\pi_1(V)\cap\pi_1(\sigma(W))^v)^{g''}$ since $v\in\pi_1(V)$. The subgroup $\pi_1(V)\cap\pi_1(\sigma(W))^v=\pi_1(\sigma(W))$ since $\pi_1(\sigma(W))=\pi_1(W)\unlhd\pi_1(V)$, so $\pi_1(V)\cap\pi_1(\sigma(W))^g=\pi_1(V)\cap\pi_1(\sigma(W))^{g''}$ which is contained in the $\pi_1(U)$-conjugacy class represented by a labeled graph in $S$ by the induction hypothesis.

Case 4: Now assume that $g=g'g_m$ and the decomposition of $g'$ begins with some nontrivial $g_1\in\pi_1(U)$. Then $\pi_1(U)\cap\pi_1(\sigma(W))^g=(\pi_1(U)\cap\pi_1(\sigma(W))^{g'})^{g_m}$ since $g_m\in\pi_1(U)$. By the induction hypothesis, $\pi_1(U)\cap\pi_1(\sigma(W))^{g'}$ is contained in a $\pi_1(U)$-conjugacy class represented by a labeled graph $H^4_U\in S$. Therefore $\pi_1(U)\cap\pi_1(W)^g$ is also represented by $H^4_U\in S$.

Since $g'$ begins with the nontrivial $g_1\in\pi_1(U)$, the path in $T$ stabilized by $\pi_1(V)\cap\pi_1(\sigma(W))^g$ passes through the vertex in $T$ stabilized by $\pi_1(U)$. Therefore $\pi_1(V)\cap\pi_1(\sigma(W))^g=\pi_1(V)\cap\pi_1(U)\cap\pi_1(\sigma(W))^g$. The previous step tells us that $\pi_1(U)\cap\pi_1(\sigma(W))^g$ is represented by the labeled graph $H^4_U\in S$, making $\pi_1(U)\cap\pi_1(\sigma(W))^g\leq\pi_1(U)$. Since $\pi_1(U)\cap\pi_1(V)=\pi_1(\sigma(W))$, $\pi_1(V)\cap\pi_1(\sigma(W))^g=\pi_1(\sigma(W))\cap\pi_1(\sigma(W))^g$, forcing $\pi_1(V)\cap\pi_1(\sigma(W))^g$ to be contained in the $\pi_1(U)$-conjugacy class represented by a connected component of $\sigma(W)\otimes_UH_U^4\in S$.

Therefore $\pi_1(U)\cap\pi_1(\sigma(W))^g$ and $\pi_1(V)\cap\pi_1(\sigma(W))^g$ are each contained in a $\pi_1(U)$-conjugacy class represented by a labeled graph in $S$ for all $g\in A_{2,2,3n}$.
\end{proof}

\begin{lemma}
    If $S$ is finite, then $A_{2,3,2n}$ has finite stature with respect to $\{\pi_1(U),\pi_1(V)\}$.
\end{lemma}
\begin{proof}
    By Proposition 3.5 in \cite{jankiewicz2023finite}, it suffices to show that there are finitely many $A_{2,3,2n}$-conjugacy classes of stabilizers of finite paths in the Bass-Serre tree $T$ that contain the vertices in $T$ stabilized by $\pi_1(U)$ and $\pi_1(V)$ respectively. Every such path stabilizer is either of the form $\pi_1(U)$, $\pi_1(V)$, $\pi_1(U)\cap\pi_1(W)^{g_1}$, $\pi_1(V)\cap\pi_1(W)^{g_2}$ or an intersection of stabilizers of these forms, for some $g_1,g_2\in A_{2,3,2n}$. By Lemma 3.8, we know that every such path stabilizer is contained in a $\pi_1(U)$-conjugacy class represented by a graph in $S$. 

    Since $\pi_1(V)/\pi_1(W)=\{[1],[v]\}$, this makes $\pi_1(V)=\pi_1(W)\bigsqcup v\pi_1(W)$. The $\pi_1(W)$-conjugacy classes of path stabilizers in $T$ are represented by the same labeled graphs in $S$ as the $\pi_1(U)$-conjugacy classes since $\pi_1(W)=\pi_1(\sigma(W))<\pi_1(U)$. The $v\pi_1(W)$-conjugacy classes are represented by the image under $\beta$ of the labeled graphs in $S$ by the definition of $\beta$. Therefore $S$ being finite, closed under fiber products, and closed under $\beta$ guarantees that $A_{2,3,2n}$ has finite stature with respect to $\{\pi_1(U),\pi_1(V)\}$.
\end{proof}

The general procedure for constructing $S$ is as follows. Begin with $S=\{\sigma(W)\}$, then:
\begin{itemize}
    \item Perform $H\otimes_{U}H$, $H\otimes_{U}K$ and $K\otimes_{U}H$ for each $H,K\in S$. Replace each resulting vertex-tuple with its second component and add to $S$ any resulting connected component that is not a labeled subgraph of an element already in $S$.
    \item Calculate $\beta(H)$ for every $H\in S$. Add $\beta(H)$ to $S$ if it is not a labeled subgraph of any element in $S$.
    \item Repeat the above two steps until no new graphs can be added to $S$ in this fashion.
\end{itemize}

\subsection{$q$-connectedness}
We introduce in this subsection the notion of $q$-connectedness, which will greatly aid in our proof that $S$ is finite in Sections 4 and 5. For the remainder of the paper, the word ``loop" will refer to a reduced loop, meaning that we remove any instances of a path going back and forth across the same edge, unless explicitly specified. We begin by defining the concept of $\beta$-consistent paths which will be the paths relevant to the remainder of the subsection.

\begin{definition}
    A path $p$ in a labeled graph $P\looparrowright\sigma(W)$ is $\beta$-consistent if $p$ is one of the following paths:
    \begin{itemize}
        \item $n$ consecutive red edges,
        \item three consecutive green edges,
        \item a green edge followed by a yellow edge entering the same vertex,
        \item a green edge followed by a yellow edge exiting the same vertex.
    \end{itemize}
\end{definition}

Figure 9 shows that the image under $\beta$ of a $\beta$-consistent loop remains $\beta$-consistent. The calculations in Figure 9 can be checked by referring to Figure 8. Notice that we omit any edges in $\sigma(W)$ that do not contribute to loops. Figure 9 shows that if we are presented with a labeled graph $Z\looparrowright\sigma(W)$ in which every simple loop is $\beta$-consistent, then every simple loop in $\beta(Z)$ will be $\beta$-consistent. 

\begin{figure}[h!]
    \centering
    \includegraphics[scale=0.3]{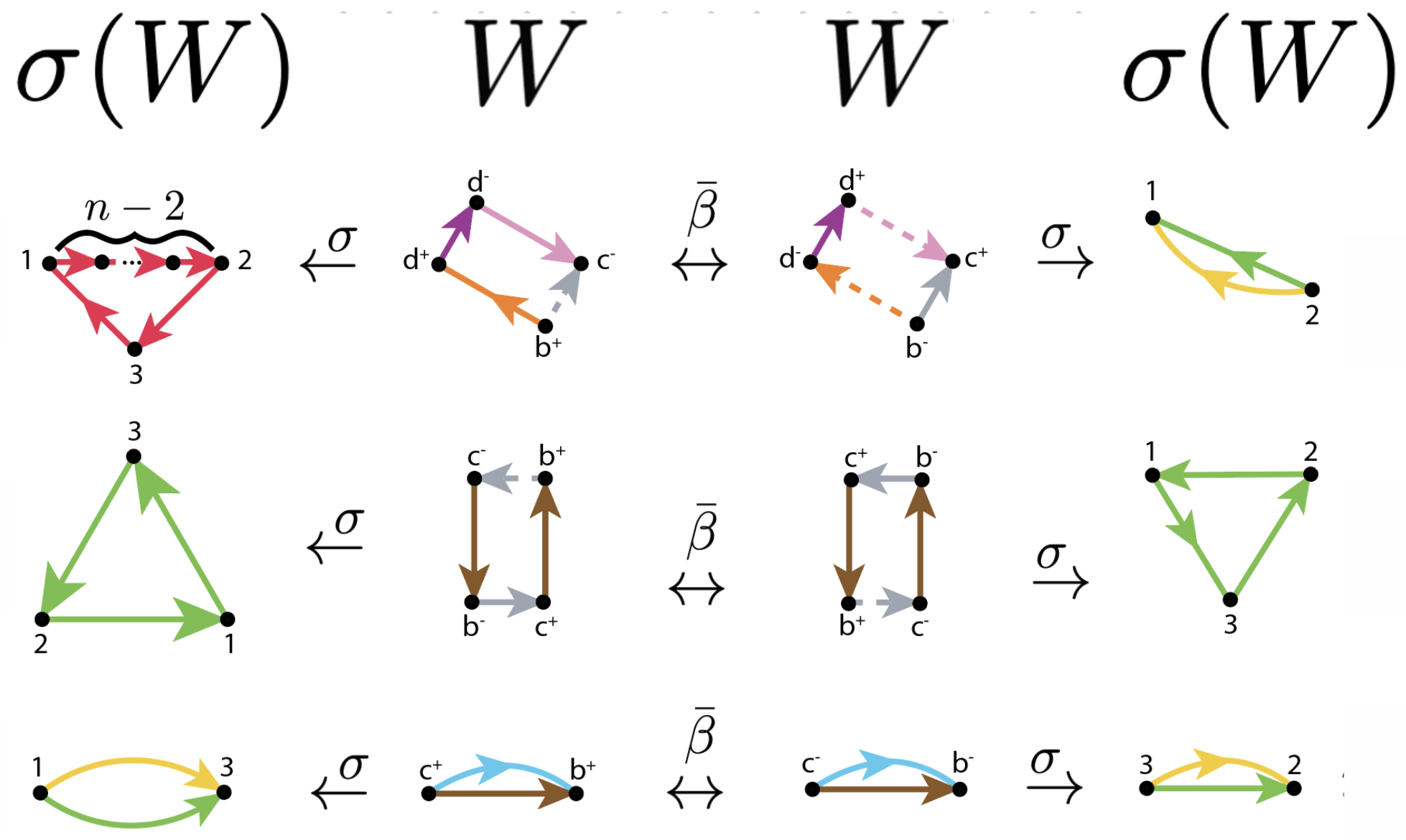}
    \caption{$\beta$ applied to $\beta$-consistent loops.}
    \label{fig:enter-label}
\end{figure}

\begin{definition}
    A connected labeled graph $K\looparrowright \sigma(W)$ is \emph{$q$-fillable} if every $\beta$-consistent path in the core of $K$ is a loop. When a labeled graph $K$ is $q$-fillable, define $q(K)$ to be the 2-complex with $K$ as its 1-skeleton and a 2-cell filling every red $n$-gon, green triangle and yellow-green bigon.
\end{definition}

For example, $\sigma(W)$ is $q$-fillable since every $\beta$-consistent path that is present in $\sigma(W)$ is closed. Notice that if we extend $\sigma(W)$ to a cover of $U$ then this extension is no longer $q$-fillable due to the non closed length-3 green paths extending from the trees attached to the vertices in the red $n-2$ path in $\sigma(W)$, for example. It is important to remember that $q$-fillability is defined with respect to the core of a labeled graph.

\begin{lemma}
    Suppose that $H,K\looparrowright \sigma(W)$ are $q$-fillable. Then:
    \begin{itemize}
        \item every connected component of $H\otimes_{U}K$ is $q$-fillable,
        \item $\beta(H)$ is $q$-fillable.
    \end{itemize}
\end{lemma}
\begin{proof}
Let $L$ be a connected component of $H\otimes_{U}K$. Then $L$ combinatorially immerses into both $H$ and $K$. Let $\ell$ be a path based at $(h,k)$ in $L$ whose edge labels form a $\beta$-consistent path. The fiber product definition allows us to view $\ell$ also as a path based at $h$ in $H$ and a path based at $k$ in $K$. Since $H$ and $K$ are $q$-fillable, this makes $\ell$ a loop in both $H$ and $K$. Therefore $\ell$ is a loop in $L$ and $H\otimes_{U}K$ is $q$-fillable by definition.

Let $\ell$ be a path in $\beta(H)$ whose edge labels form a $\beta$-consistent path. Figures 8 \& 9 show that $\beta$ maps every $\beta$-consistent path to another $\beta$-consistent path. Therefore the preimage of $\ell$ under $\beta$, $\beta^{-1}(\ell)$, is a collection of $\beta$-consistent paths in $H$. Since $H$ is $q$-fillable, this makes every path in $\beta^{-1}(\ell)$ a loop. Therefore $\ell$ itself is also forced to be a loop, making $\beta(H)$ $q$-fillable.
\end{proof}

Notice that the construction of $S$ always begins with $\sigma(W)\otimes_U\sigma(W)$. Since $\sigma(W)$ is $q$-fillable, by Lemma 3.12 this makes every labeled graph in $S$ $q$-fillable.

\begin{definition}
    A $q$-fillable graph $K\looparrowright W$ is \emph{$q$-connected} when $q(K)$ is simply connected.
\end{definition}

\begin{remark}
    Note that a graph $K$ being $q$-connected forces every loop $\ell$ in $K$ to be a concatenation of $\beta$-consistent loops. Furthermore, if $\ell$ is a simple loop, meaning that is contains no repeated vertices or edges, then $\ell$ is $\beta$-consistent.
\end{remark}

There will be many examples of $q$-connected graphs and graphs that are not $q$-connected in Sections 4 \& 5. An example of a graph that is not $q$-connected is $\sigma(W)$ due to the presence of homotopically nontrivial red-green bigons. The presence of a red-green bigon is often what keeps graphs from being $q$-connected, though this is not always the case. A graph being $q$-connected means that it is essentially a ``dead end" with regards to further $\beta$ and fiber product calculations. We use the remainder of this section to prove that this extremely helpful property holds.

\begin{lemma}
    Let $H,K\looparrowright \sigma(W)$ be $q$-fillable graphs where $K$ is $q$-connected. Then every connected component of $K\otimes_{U} H$ and $H\otimes_{U} K$ is a subgraph of $K$.
\end{lemma}
\begin{proof}
Let $L$ be a connected component of $K\otimes_{U}H$. Then there is a combinatorial immersion $\phi:L\to K$. In order for $\phi$ to not be an embedding, there needs to exist a non closed path $p$ in $L$ such that $\phi(p)$ is a loop in $K$. Assume by way of contradiction that such a non closed path $p$ exists. Since $K$ is $q$-connected, this makes $\phi(p)$ a loop that arises as a concatenation of $\beta$-consistent loops. Therefore $p$ is a path in $L$ that arises as a concatenation of $\beta$-consistent paths. Lemma 3.12 guarantees that $L$ is $q$-fillable, forcing every path in $L$ that arises as a concatenation of $\beta$-consistent paths to be a concatenation of loops. Therefore $p$ is a loop.
\end{proof}

\begin{lemma}
If $K$ is $q$-connected, then $\beta(K)$ is $q$-connected as well.
\end{lemma}
\begin{proof}
    The labeled graph $K$ being $q$-connected means that every simple loop in $K$ is a red cycle of length $n$, a green triangle or a yellow-green bigon. As shown in Figure 9, this collection of loops is closed under $\beta$. Therefore the simple loops in $\beta(K)$ also arise as red cycles of length $n$, green triangles and yellow-green bigons, making $\beta(K)$ $q$-connected as well.
\end{proof}

\section{Residual finiteness of $A_{2,3,2n}$ for $n>4$}
The goal of this section is to prove the following theorem:
\begin{theorem}
    $A_{2,3,2n}$ for $n>4$ is residually finite.
\end{theorem}

The $n=4$ case will be proven in Section 5. Combining Lemma 3.9 with Theorem 1.3 reduces the proof of this theorem to proving that a finite set $S$, as described in Definition 3.7, exists.

\subsection{Iterative construction of $S$} 

We begin with $S=\{\sigma(W)\}$ since $S$ must, at minimum, contain $\sigma(W)$. Conjugating $\sigma(W)$ by any element of $\pi_1(U)$ does not change the structure or labels of $\sigma(W)$, it just shifts the basepoint. Also, since $\pi_1(\sigma(W))= \pi_1(W)$ and $\pi_1(W)\unlhd\pi_1(V)$, $v^{-1}\pi_1(\sigma(W))v=v^{-1}\pi_1(W)v=  \pi_1(W)=\pi_1(\sigma(W))$ where $v$ is the fixed nontrivial coset representative of $\pi_1(V)/\pi_1(W)\cong\mathbb{Z}/2\mathbb{Z}$. Pictorially, $\beta(\sigma(W))\cong \sigma(W)$ is shown in Figure 10.
\begin{figure}[h]
    \centering
    \includegraphics[scale=0.35]{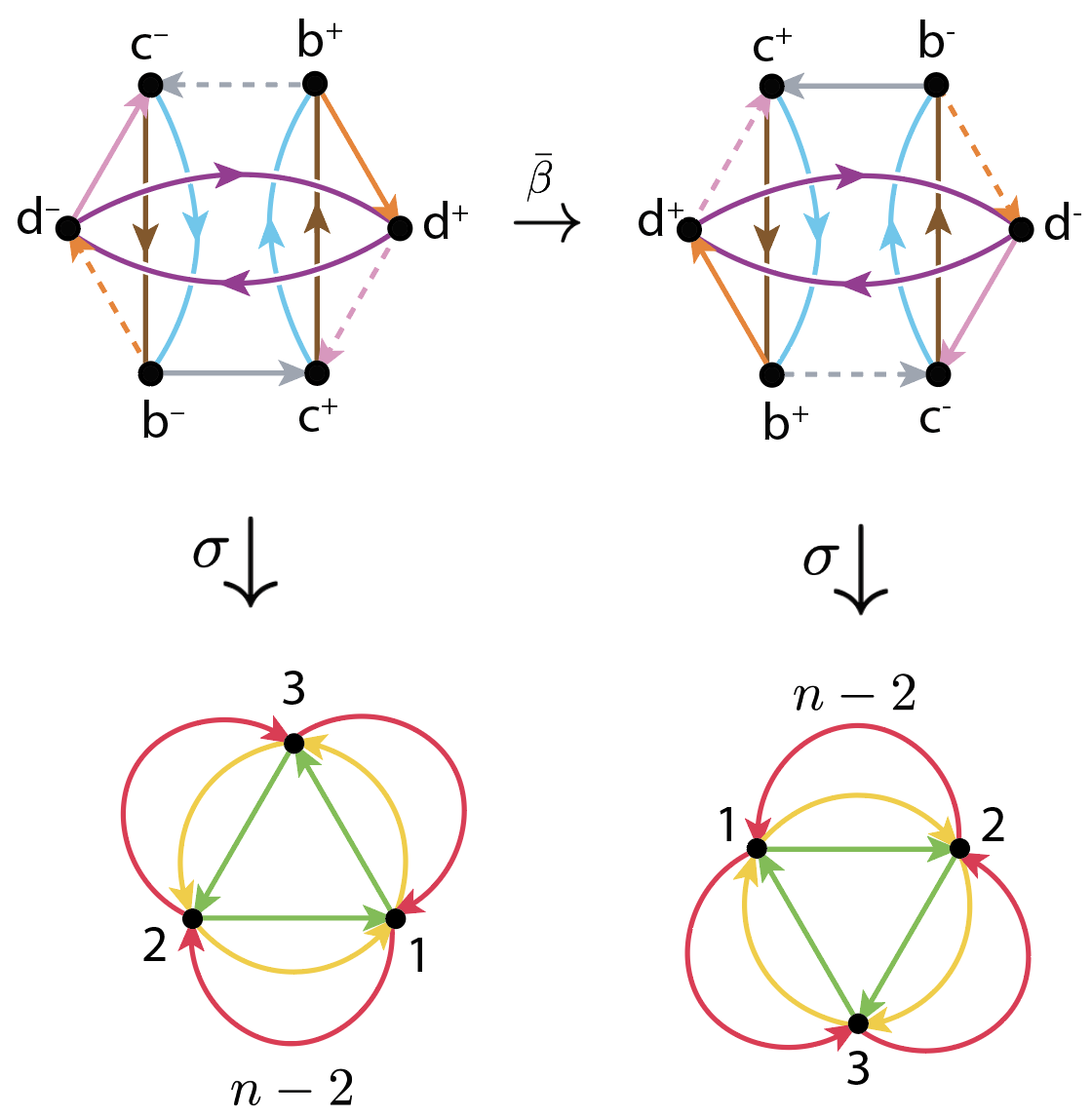}
    \caption{$\beta(\sigma(W))$ is a rotated copy of $\sigma(W)$}
    \label{fig:enter-label}
\end{figure}

We end this subsection with the statement of the lemma that we will spend the rest of the section proving.

\begin{lemma}
The set $S$ is finite for $A_{2,3,2n}$ with $n>4$.
\end{lemma}

So far $S=\{\sigma(W)\}$ is closed under basepoint translation, but we also $S$ need to be closed under fiber product. 

\subsection{$\sigma(W)\otimes_{U}\sigma(W)$ computation}

\begin{lemma}
    For $A_{2,3,2n}$ with $n>4$, $\sigma(W)\otimes_{U}\sigma(W)$ has the following connected components:\begin{itemize}
        \item one copy of $\sigma(W)$,
        \item $n-5$ copies of an $n$-gon with red edges,
        \item one of each of the graphs in Figure 11, to be denoted henceforth by $X_1$ (leftmost) and $X_2$ (rightmost).
    \end{itemize}
\end{lemma}
\begin{figure}[h]
    \centering
    \includegraphics[scale=0.3]{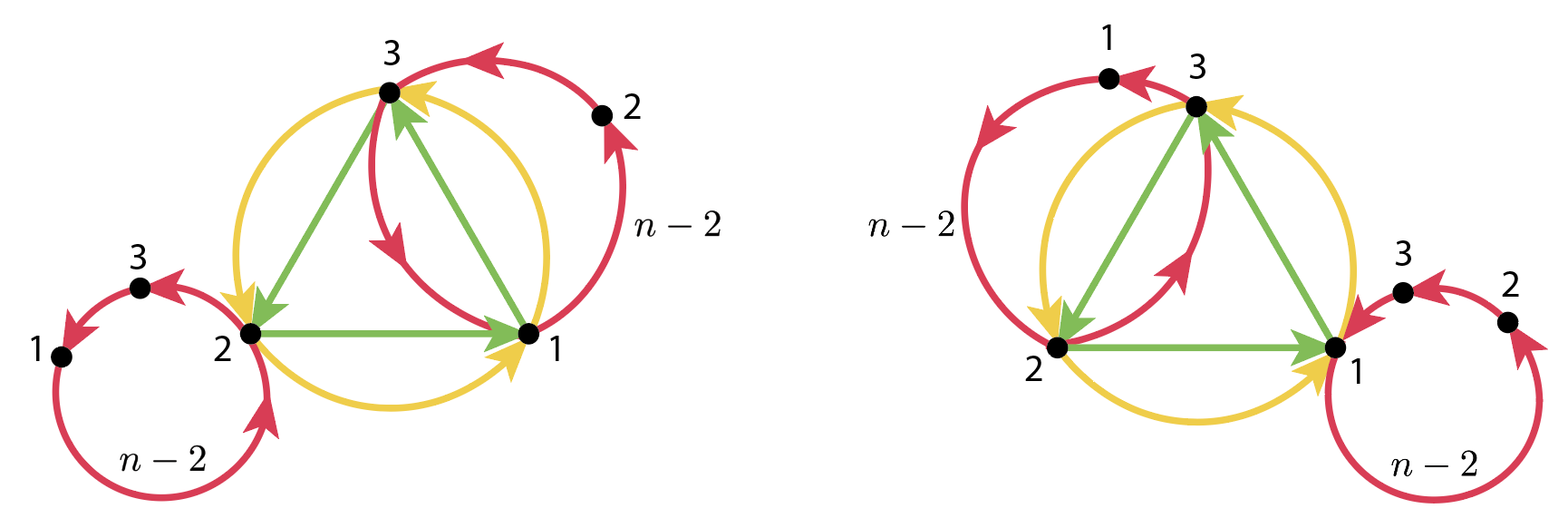}
    \caption{$X_1$ \& $X_2$}
    \label{fig:enter-label}
\end{figure}

\begin{proof}
    We start by focusing on the edges found in $\sigma(W)\otimes_{U}\sigma(W)$. The labeled graph $\sigma(W)$ has three yellow edges, three green edges and $n$ red edges. Therefore $\sigma(W)\otimes_{U}\sigma(W)$ must have nine yellow edges, nine green edges and $n^2$ red edges. By direct computation, these yellow edges and green edges arise as the triangles 
    \begin{itemize}
        \item $(1,1)\to(3,3)\to(2,2)\to(1,1)$
        \item $(1,2)\to(3,1)\to(2,3)\to(1,2)$
        \item $(1,3)\to(3,2)\to(2,1)\to(1,3)$
    \end{itemize}
    
    The first triangle is present in the connected component that is identical to $\sigma(W)$. In this connected component every vertex is labeled with a tuple whose elements are identical. After restricting each vertex label to the second component in its tuple, the other two triangles become the triangles present in $X_1$ and $X_2$ respectively. The placement of the red edges in $\sigma(W)$, $X_1$ and $X_2$ follow directly from these vertex-labelings. 
    
    We now shift our focus to the vertices present in $\sigma(W)\otimes_{U}\sigma(W)$. Since $U$ has only one vertex and the labeled graph $\sigma(W)$ has $n$ vertices, $\sigma(W)\otimes_{U}\sigma(W)$ has $n^2$ vertices. The copy of $\sigma(W)$ that arises as a connected component of $\sigma(W)\otimes_{U}\sigma(W)$ has $n$ vertices, $X_1$ has $2n$ vertices and $X_2$ has $2n$ vertices. Therefore there are $n^2-5n$ vertices unaccounted for. Since we have already analyzed all of the yellow and green edges in $\sigma(W)\otimes_{U}\sigma(W)$, only red edges can connect the remaining $n^2-5n$ vertices. Every vertex in $\sigma(W)$ is present in a red $n$-gon. Therefore every vertex in $\sigma(W)\otimes_{U}\sigma(W)$ must also be present in a red $n$-gon. Red $n$-gons have $n$ edges, therefore the remaining $n^2-5n$ vertices must be distributed across $n-5$ red $n$-gons.
\end{proof}

\subsection{The proof of Lemma 4.2}
\begin{proof}
Lemma 4.3 describes the connected components of $\sigma(W)\otimes_{U}\sigma(W)$ as being a copy of $\sigma(W)$, $X_1$ and $X_2$, along with a collection of red $n$-gons. Let $R$ be a labeled graph that is just a red $n$-gon and $K\looparrowright U$ a labeled graph. Since $R$ is a subgraph of $\sigma(W)$, by Lemma 3.3, every connected component of $R\otimes_{U}K$ is a subgraph of $\sigma(W)\otimes_{U}K$, so we do not add these components to $S$.

So far $S=\{\sigma(W),X_1,X_2\}$. We need $S$ to be closed under $\beta$. Figure 10 shows that $\beta(\sigma(W))=\sigma(W)$, so it remains to calculate $\beta(X_1)$ and $\beta(X_2)$. These calculations are carried out in Figures 12 \& 13 by applying $\bar{\beta}$ and $\sigma$ to each edge. Refer to Figure 8 for more details about how to calculate the image of each edge. These calculations result in two new labeled graphs that must be included in $S$. The vertices $\bar{1}$ (resp. $\bar{2}$) in $\beta(X_1)$ and $\beta(X_2)$ is a preimage of the vertex labeled 1 (resp. 2) in $\sigma(W)$, and the notation is used to differentiate the two preimages for ease of computation later on.

\begin{figure}[h]
    \centering
    \includegraphics[scale=0.3]{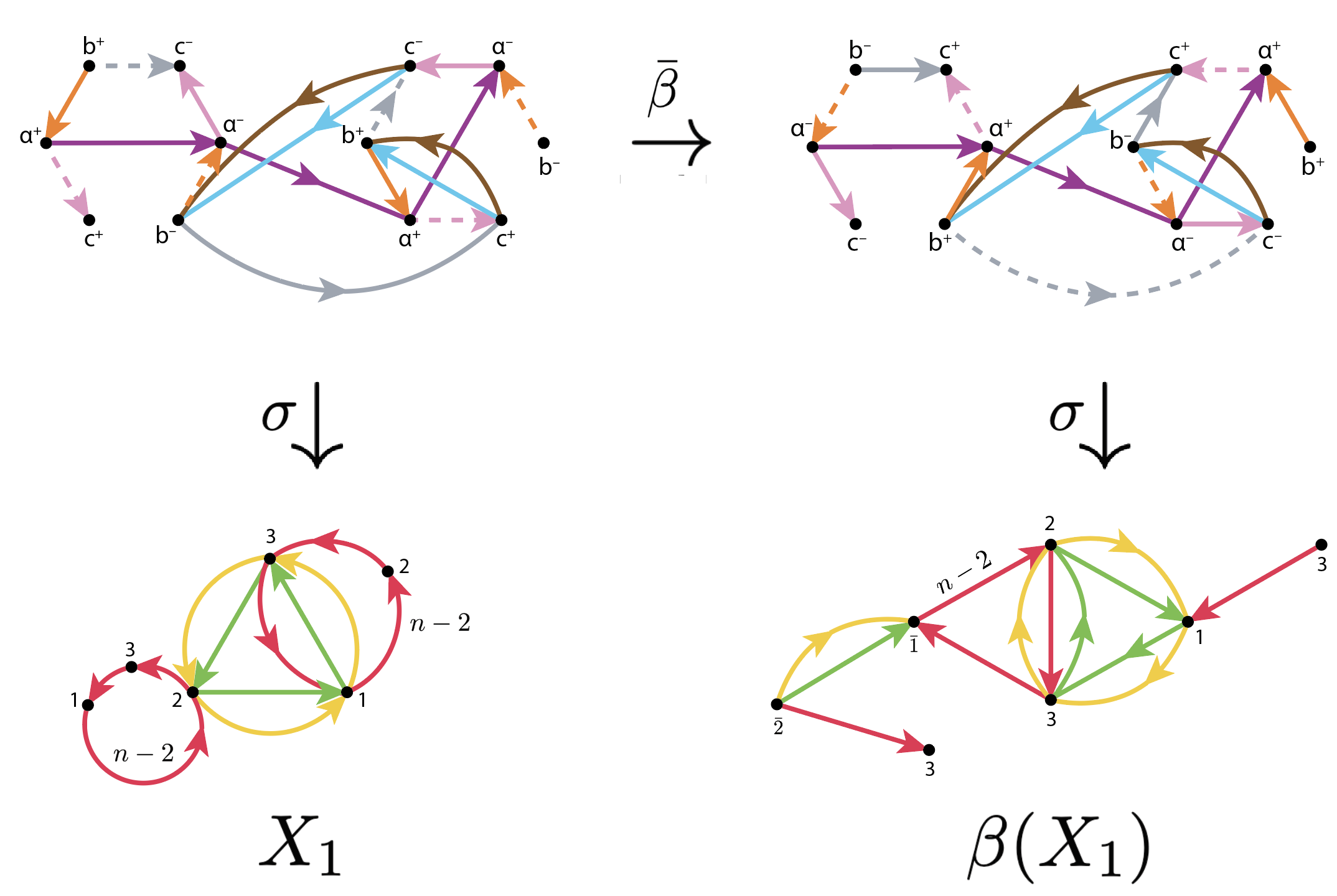}
    \caption{$\beta(X_1)$}
    \label{fig:enter-label}
\end{figure}
\begin{figure}[h!]
    \centering
    \includegraphics[scale=0.35]{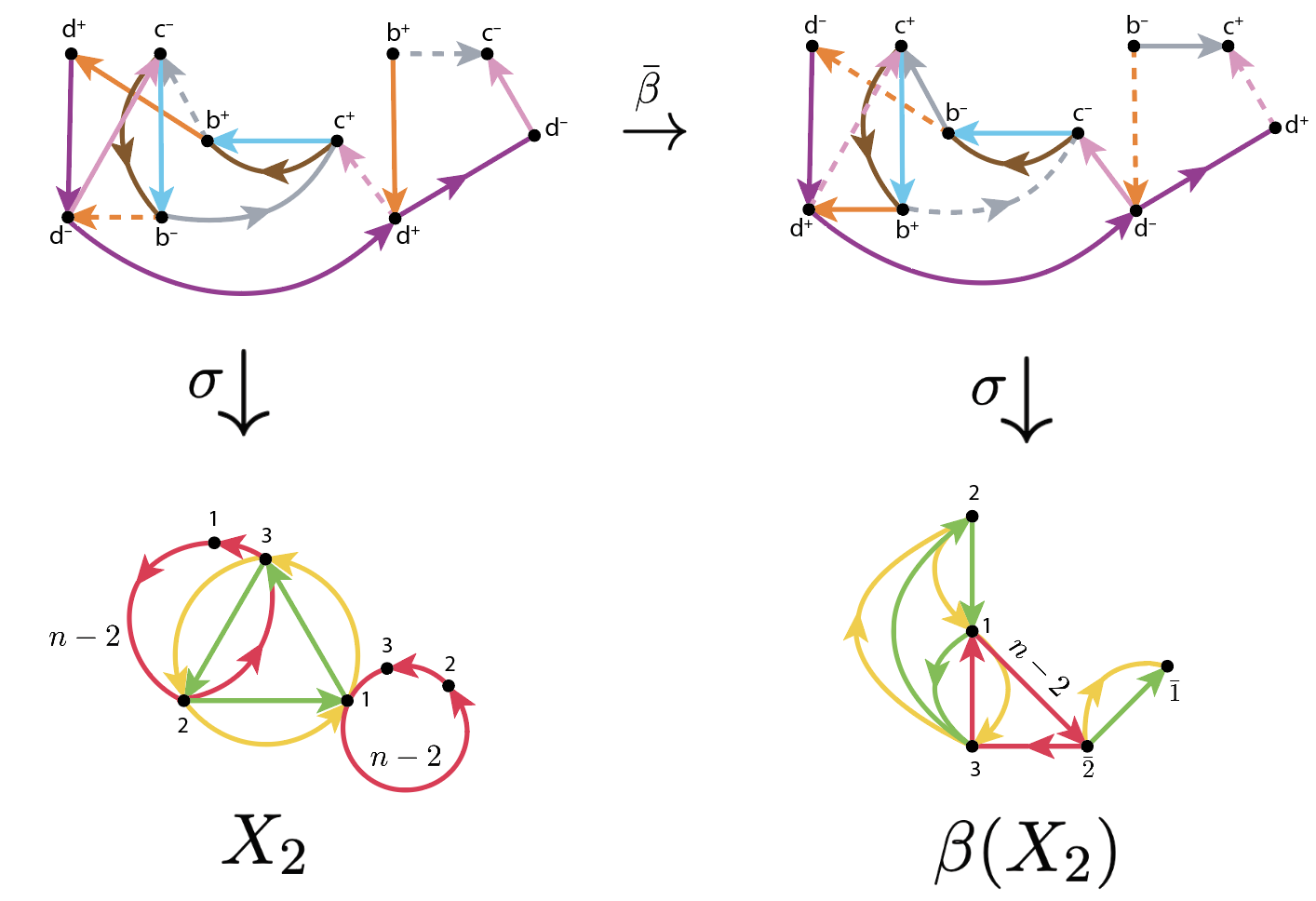}
    \caption{$\beta(X_2)$}
    \label{fig:enter-label}
\end{figure}

We now have $S=\{\sigma(W),X_1,X_2,\beta(X_1),\beta(X_2)\}$, a set closed under $\beta$. Lemmas 3.14 \& 3.15 will be the tools we use to work with $q$-connected components separately from those that are not $q$-connected. All $q$-connected components that arise in the following calculations will be approached in detail in Section 4.3.5. In Sections 4.3.1-4.3.4 we focus on showing that every connected component of each pairwise fiber product of elements currently in $S$ is either a subgraph of a graph in $S$ or is $q$-connected.  To prove this claim we proceed in a similar fashion as in Lemma 4.3. As $n$ increases, the number of red edges increases, but the number of yellow and green edges remains the same, so it suffices to focus on how the yellow and green edges are dispersed among the connected components. The collection of connected components of a fiber product that contain yellow and green edges will be termed the \emph{relevant portion} for the remainder of the proof. The connected components that are not relevant therefore only consist of red edges, making them subgraphs of $\sigma(W)$ and thus unable to introduce any new graphs to $S$. 

\subsubsection{All fiber products involving only $\sigma(W)$, $X_1$ and $X_2$}
We begin by calculating the relevant portion of $\sigma(W)\otimes_{U}X_1$. By direct computation, these components are two copies of $X_1$ and a graph that is $q$-connected, as seen in Figure 14.

\begin{figure}[h!]
            \centering
            \includegraphics[scale=0.3]{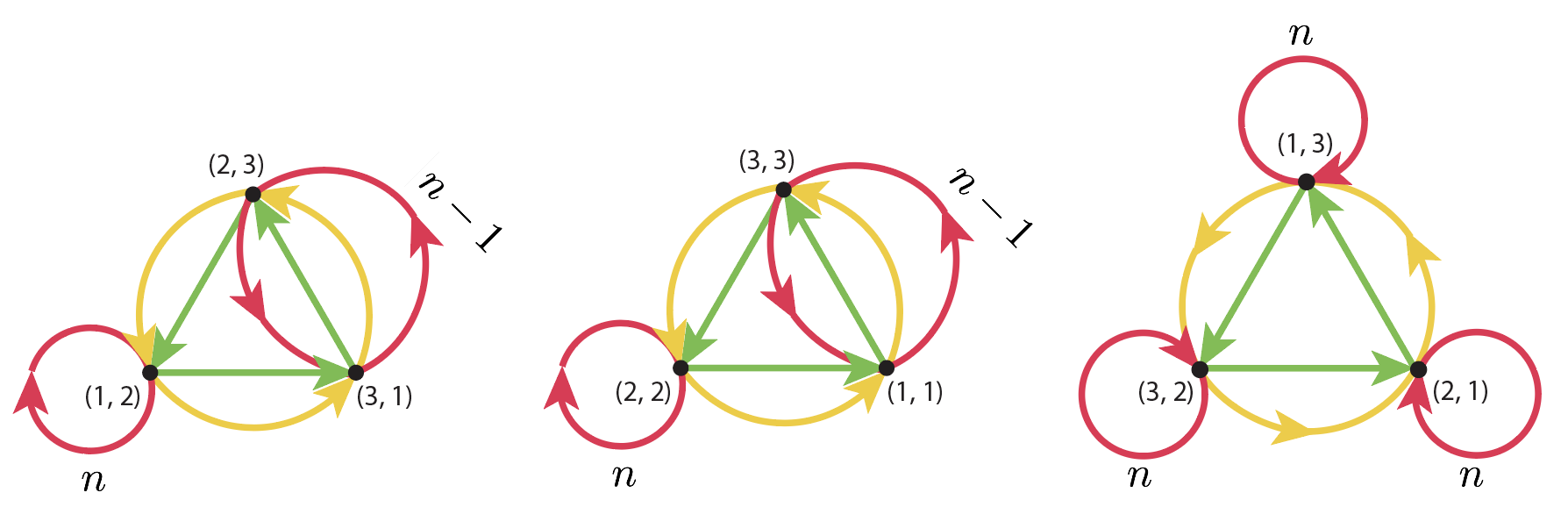}
            \caption{The relevant portion of $\sigma(W)\otimes_{U}X_1$}
            \label{fig:enter-label}
        \end{figure}

We compute the relevant portion of $X_1\otimes_{U}\sigma(W)$ by swapping the entries in each vertex-tuple in Figure 14, resulting in $X_1$, $X_2$ and a $q$-connected graph.

The fiber products used to build $S$ are with respect to $U$, meaning that we view each component of the fiber product as being combinatorially immersed in $U$. The graph $U$ has only one vertex, so the vertex labels of the graphs in $S$ are irrelevant to the fiber product calculations. However, the vertex labels are extremely important when it comes to computing the image of a graph under $\beta$. Notice that $X_1$ and $X_2$ are identical as unlabeled graphs. We denote this property by $X_1\cong X_2$. Since $X_1\cong X_2$, the connected components that arise from $Y\otimes_{U}X_1$ (resp. $X_1\otimes_{U}Y$)  will be the same unlabeled graphs as the connected components in $Y\otimes_{U}X_2$ (resp. $X_2\otimes_{U}Y$) for all $Y\looparrowright U$. 

In particular, the connected components of $\sigma(W)\otimes_{U}X_2$ are identical to the connected components of $\sigma(W)\otimes_{U}X_1$ up to vertex labels. Since a graph is $q$-connected regardless of its vertex labels, it suffices to calculate the vertex labels of the two non $q$-connected components of $\sigma(W)\otimes_{U}X_2$. This results in two copies of $X_2$ as seen in Figure 15.

\begin{figure}[h!]
        \centering
        \includegraphics[width=0.5\linewidth]{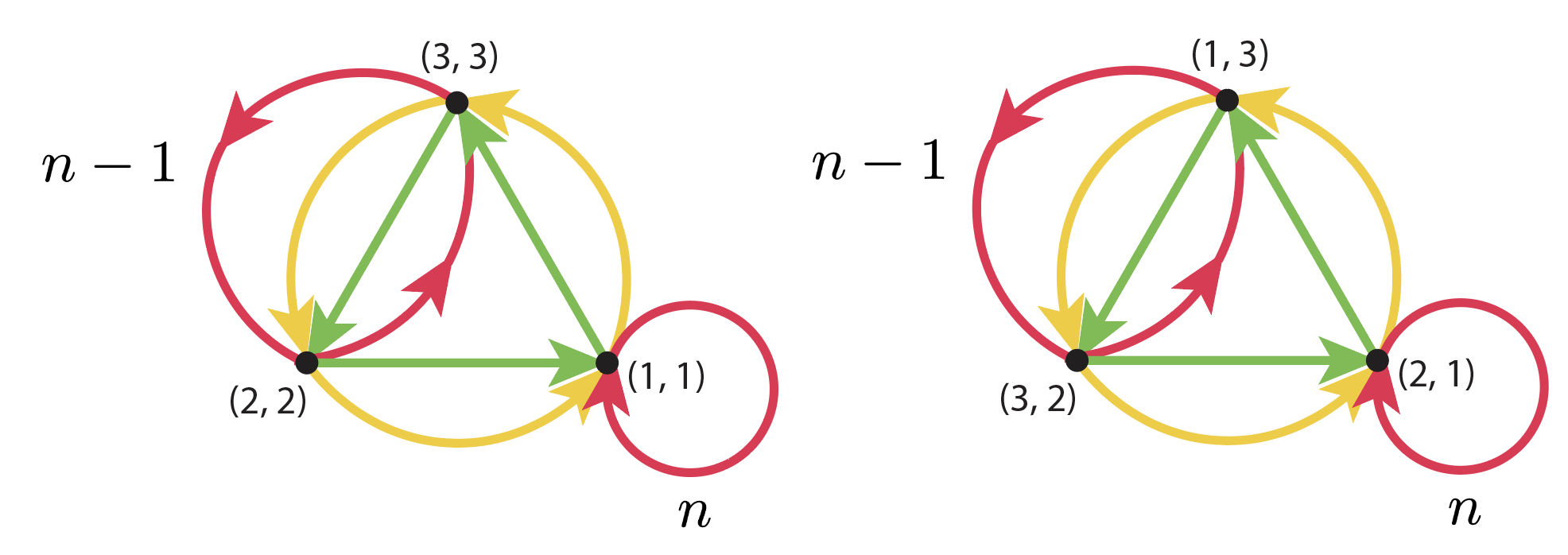}
        \caption{The non $q$-connected relevant portion of $\sigma(W)\otimes_{U}X_2$}
        \label{fig:enter-label}
    \end{figure} 
To calculate the vertex labels of the non $q$-connected relevant portion of $X_2\otimes_{U}\sigma(W)$, we swap the tuple-entries at each vertex in Figure 15, resulting in a copy of $X_2$ and $X_1$.

The relevant portion of $X_1\otimes_{U} X_2$ is 2 $q$-connected graphs and $X_2$, as seen in Figure 16.
\begin{figure}[h!]
            \centering
            \includegraphics[scale=0.25]{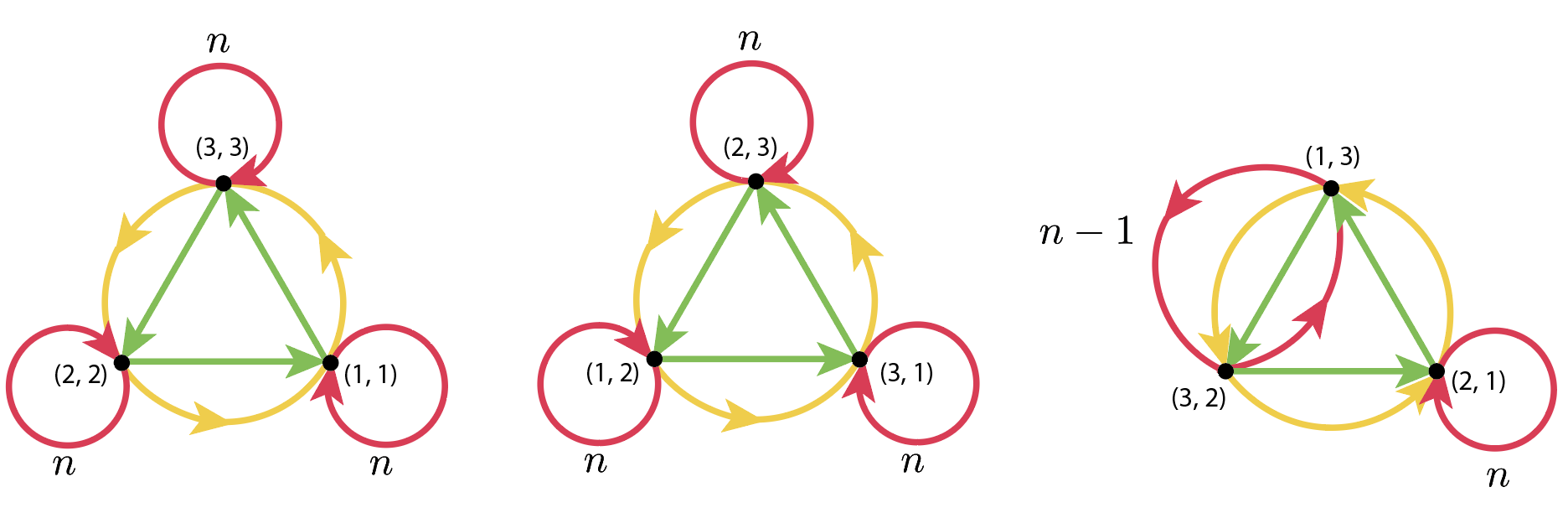}
            \caption{The relevant portion of $X_1\otimes_{U}X_2$}
            \label{fig:enter-label}
        \end{figure}

The rightmost graph in Figure 16 is the only graph in the relevant portion of $X_1\otimes_{U}X_2$ to not be $q$-connected due to the presence of the red-green bigon between the vertices labeled $(1,3)$ and $(3,2)$. To obtain the non $q$-connected relevant portion of $X_2\otimes_{U}X_1$, we swap the tuple-entries at every vertex in the rightmost graph in Figure 16, resulting in a copy of $X_1$.

Since $X_1\cong X_2$, both $X_1\otimes_{U}X_1$ and $X_2\otimes_{U}X_2$ have the same number of non $q$-connected relevant components as $X_1\otimes_{U}X_2$, which has only one. For $X_1\otimes_{U}X_1$, that connected component is the copy of $X_1$ shown in Figure 17.

\begin{figure}[h!]
        \centering
        \includegraphics[scale=0.2]{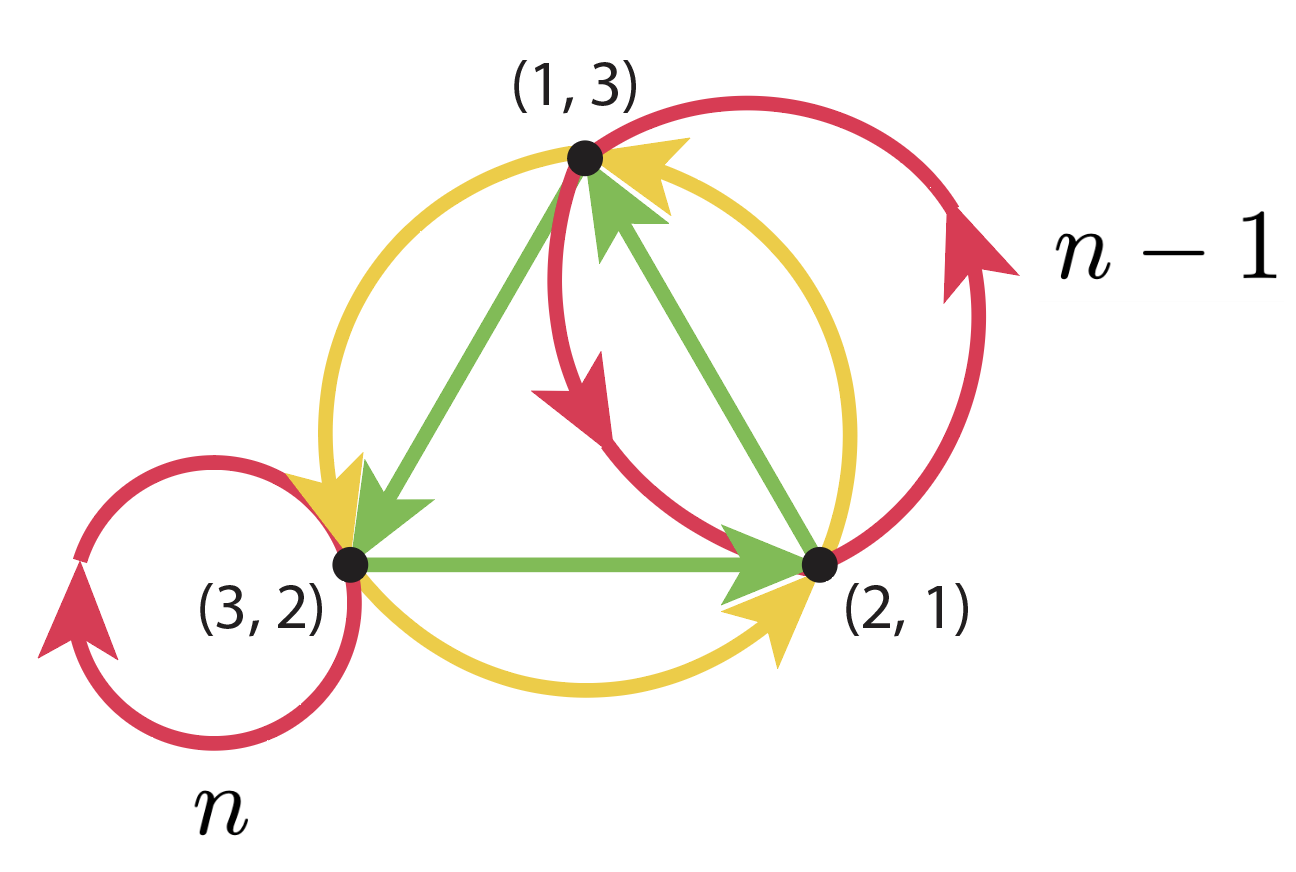}
        \caption{The non $q$-connected relevant portion of $X_1\otimes_{U}X_1$}
        \label{fig:enter-label}
    \end{figure}

For $X_2\otimes_{U}X_2$, the non $q$-connected relevant connected component is the copy of $X_2$ shown in Figure 18.

\begin{figure}[h!]
        \centering
        \includegraphics[scale=0.25]{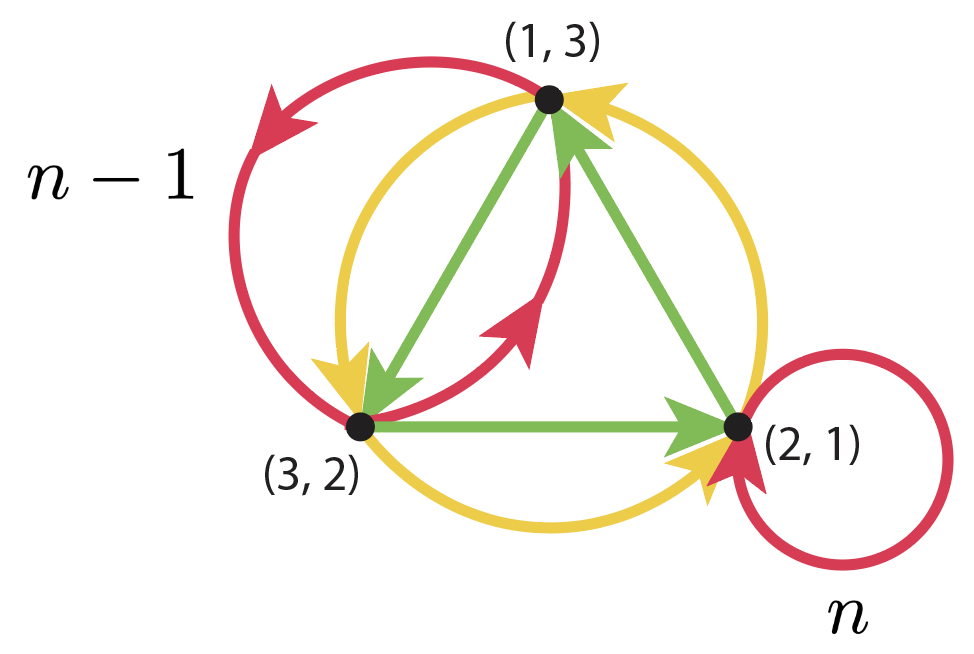}
        \caption{The non $q$-connected relevant portion of $X_2\otimes_{U}X_2$}
        \label{fig:enter-label}
    \end{figure}

\subsubsection{$Z\otimes_{U}\beta(X_1)$ where $Z\in\{\sigma(W),X_1,X_2\}$}

The relevant portion of $\sigma(W)\otimes_{U}\beta(X_1)$ is a copy of $\beta(X_1)$, two $q$-connected graphs and a subgraph of $X_2$, as shown in Figure 19.
     \begin{figure}[h!]
        \centering
        \includegraphics[scale=0.45]{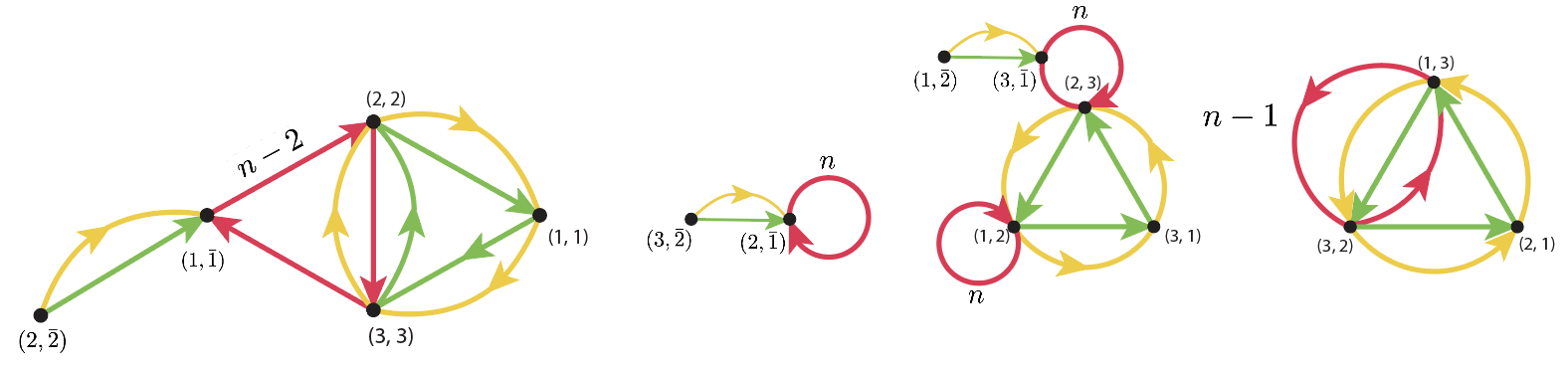}
        \caption{The relevant portion of $\sigma(W)\otimes_{U}\beta(X_1)$}
        \label{fig:enter-label}
    \end{figure}
     
Swapping the tuple-entries at each vertex in the leftmost and rightmost graphs in Figure 19 results in the non $q$-connected relevant portion of $\beta(X_1)\otimes_{U}\sigma(W)$ being $\beta(X_1)$ and a subgraph of $X_1$.

The relevant portion of $X_1\otimes_{U}\beta(X_1)$ consists of four $q$-connected components and a subgraph of $X_2$ as shown in Figure 20.
    \begin{figure}[h!]
        \centering
        \includegraphics[scale=0.4]{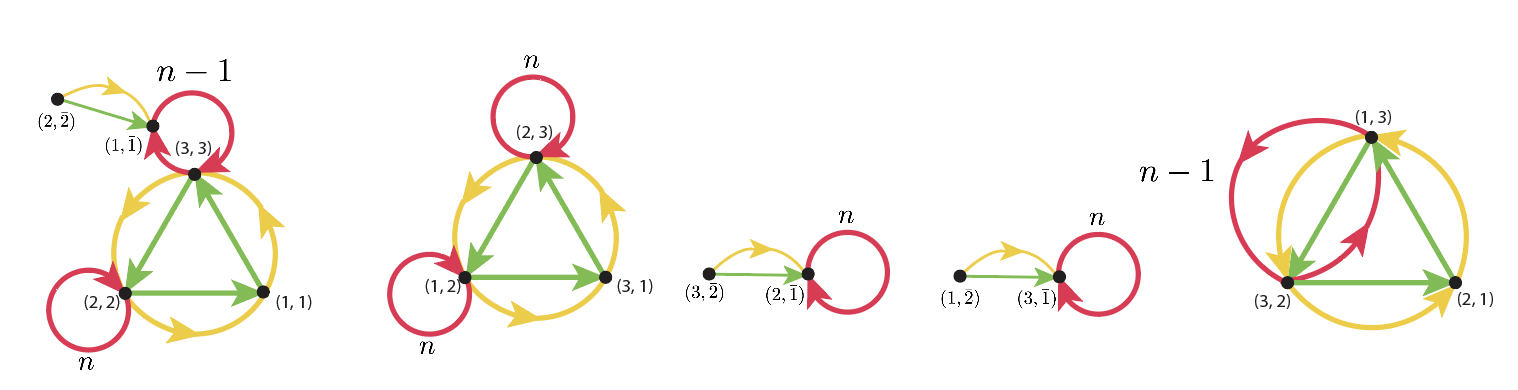}
        \caption{The relevant portion of $X_1\otimes_{U}\beta(X_1)$}
        \label{fig:enter-label}
    \end{figure} 

    Swapping the tuple-entries at each vertex in the rightmost graph in Figure 20 results in the only non $q$-connected relevant component of $\beta(X_1)\otimes_{U}X_1$ being a subgraph of $X_1$.

Since $X_1\cong X_2$, we know from the above calculations that there is only one non $q$-connected relevant component of $X_2\otimes_{U}\beta(X_1)$. This component is the subgraph of $X_2$ shown in Figure 21. Swapping the tuple-entries at each vertex leaves the labels unchanged, so the only non $q$-connected relevant component of $\beta(X_1)\otimes_{U}X_2$ is the same subgraph of $X_2$.

\begin{figure}[h!]
            \centering
            \includegraphics[scale=0.35]{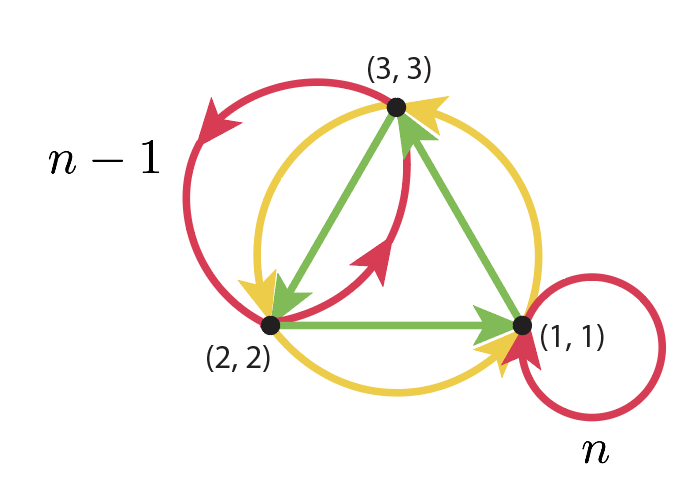}
            \caption{The relevant portion of both $X_2\otimes_{U}\beta(X_1)$ \& $\beta(X_1)\otimes_{U}X_2$}
            \label{fig:enter-label}
        \end{figure}

\subsubsection{$Z\otimes_{U}\beta(X_2)$ where $Z\in\{\sigma(W),X_1,X_2\}$}
        
The relevant portion of $\sigma(W)\otimes_{U}\beta(X_2)$ is a copy of $\beta(X_2)$, a subgraph of $X_1$ and two $q$-connected graphs as shown in Figure 22.
\begin{figure}[h!]
            \centering
            \includegraphics[scale=0.3]{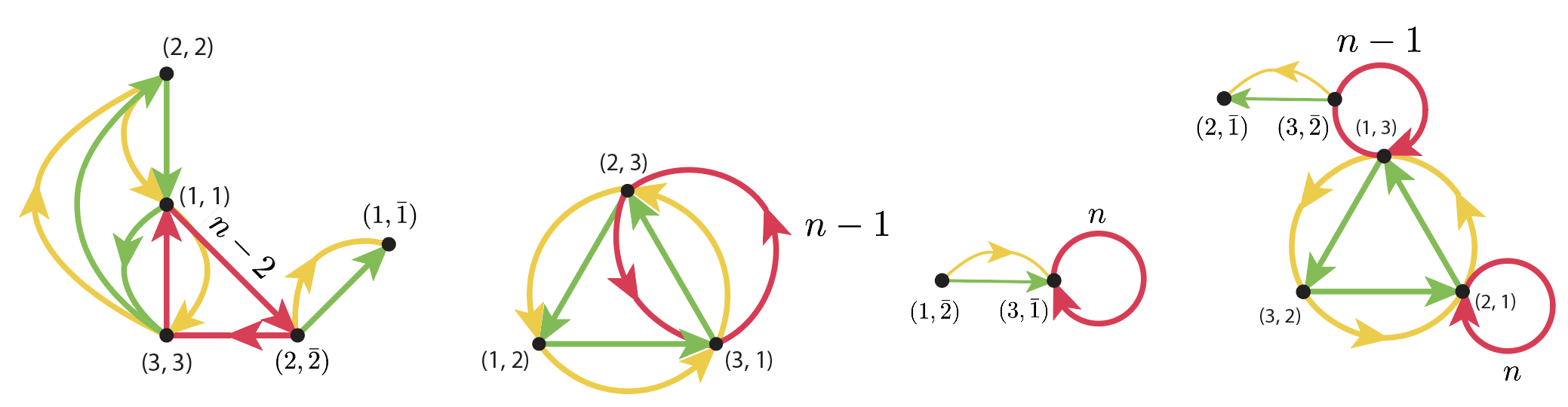}
            \caption{The relevant portion of $\sigma(W)\otimes_{U}\beta(X_2)$}
            \label{fig:enter-label}
        \end{figure}

Swapping the tuple-entries at each vertex of the two leftmost graphs in Figure 22 results in the non $q$-connected relevant portion of $\beta(X_2)\otimes_{U}\sigma(W)$ being $\beta(X_2)$ and a subgraph of $X_2$.

The relevant portion of $X_1\otimes_{U}\beta(X_2)$ consists of 4 $q$-connected graphs and a subgraph of $X_1$ as shown in Figure 23.

\begin{figure}[h!]
            \centering
            \includegraphics[scale=0.45]{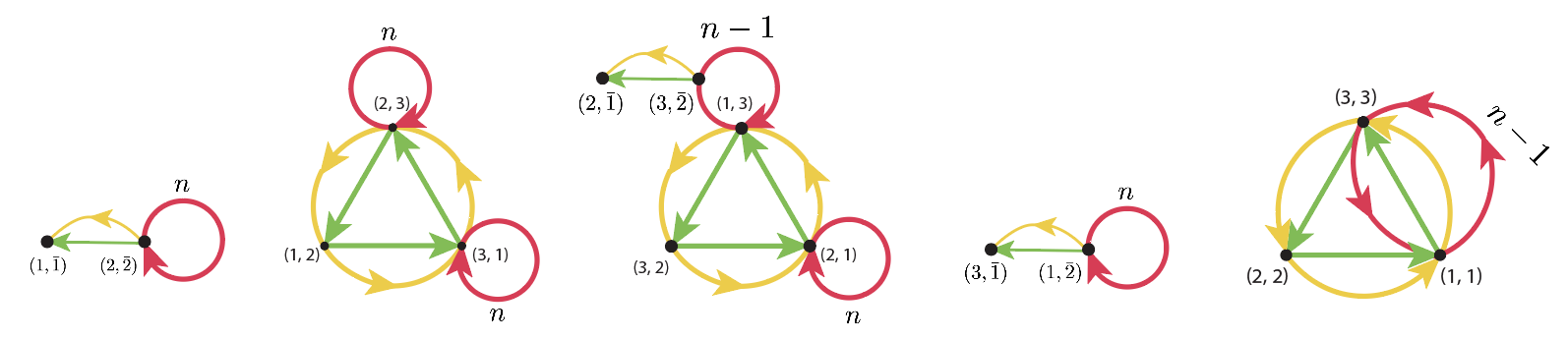}
            \caption{The relevant portion of $X_1\otimes_{U}\beta(X_2)$}
            \label{fig:enter-label}
        \end{figure}

Swapping the tuple-entries at each vertex in the rightmost graph in Figure 23 results in the same graph being the only non $q$-connected relevant component of $\beta(X_2)\otimes_{U}X_1$.

Since $X_1\cong X_2$, the only non $q$-connected relevant component of $X_2\otimes_{U}\beta(X_2)$ is the subgraph of $X_1$ shown in Figure 24.

Swapping the tuple-entries at each vertex in the graph in Figure 24 results in the only non $q$-connected relevant component of $\beta(X_2)\otimes_{U}X_2$ being a subgraph of $X_2$.

\begin{figure}[h!]
        \centering
        \includegraphics[scale=0.35]{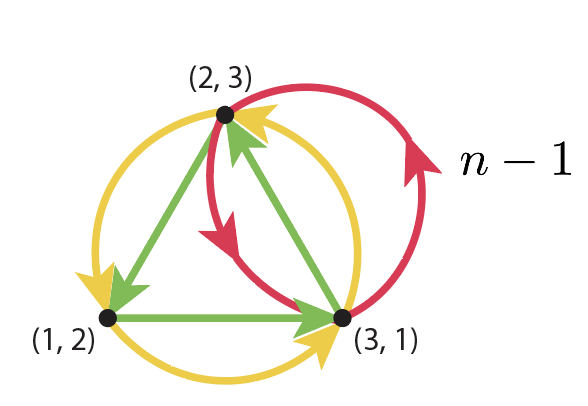}
        \caption{The relevant portion of $X_2\otimes_{U}\beta(X_2)$}
        \label{fig:enter-label}
    \end{figure}

\subsubsection{$\beta(X_1)\otimes_{U}\beta(X_1),\beta(X_1)\otimes_{U}\beta(X_2)$ and $\beta(X_2)\otimes_{U}\beta(X_2)$}

The relevant portion of $\beta(X_1)\otimes_{U}\beta(X_1)$ consists of a copy of $\beta(X_1)$ and 6 $q$-connected graphs as shown in Figure 25. 
\begin{figure}[h!]
        \centering
        \includegraphics[scale=0.4]{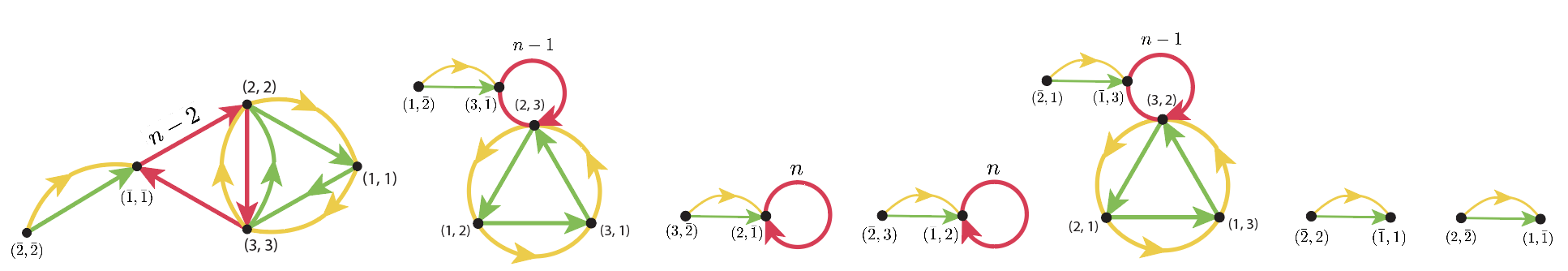}
        \caption{The relevant portion of $\beta(X_1)\otimes_{U}\beta(X_1)$}
        \label{fig:enter-label}
    \end{figure}

The relevant portion of $\beta(X_1)\otimes_{U}\beta(X_2)$ is a subgraph of $X_1$ and 6 $q$-connected graphs as shown in Figure 26. 
\begin{figure}[h!]
        \centering
        \includegraphics[scale=0.4]{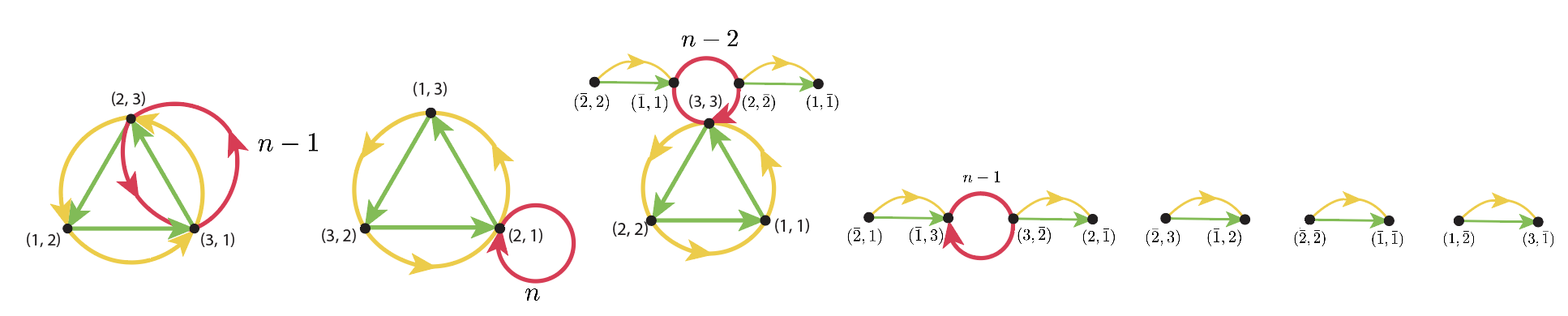}
        \caption{The relevant portion of $\beta(X_1)\otimes_{U}\beta(X_2)$}
        \label{fig:enter-label}
    \end{figure}

Swapping the tuple-entries of the leftmost graph in Figure 26 results in the only non $q$-connected relevant component of $\beta(X_2)\otimes_{U}\beta(X_1)$ being a subgraph of $X_2$. 

The relevant portion of $\beta(X_2)\otimes_{U}\beta(X_2)$ is a copy of $\beta(X_2)$ and 6 $q$-connected graphs shown in Figure 27.
\begin{figure}[h!]
        \centering
        \includegraphics[scale=0.4]{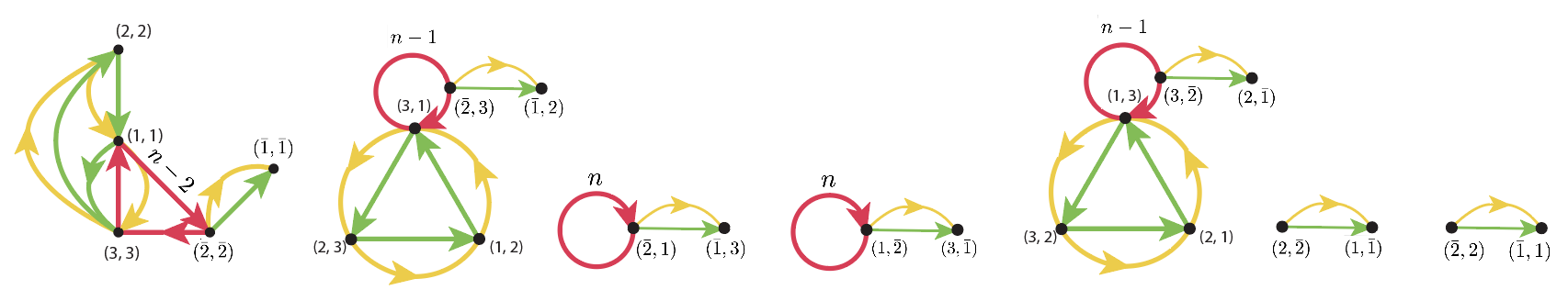}
        \caption{The relevant portion of $\beta(X_2)\otimes_{U}\beta(X_2)$}
        \label{fig:enter-label}
    \end{figure}

    \subsubsection{Assembling the pieces to finalize $S$}

    We now have a complete enumeration of every relevant connected component in every fiber product coming from pairs of graphs in $S=\{\sigma(W),X_1,X_2,\beta(X_1),\beta(X_2)\}$. Sections 4.3.1-4.3.4 show that every such component is either a subgraph of a graph in $S$ or is $q$-connected. In each fiber product computed in the last subsection, there are also non-relevant connected components that arise, consisting of unaccounted for vertices and red edges. Note that every red edge in every graph in $S=\{\sigma(W),X_1,X_2,\beta(X_1),\beta(X_2)\}$ is contained in a red $n$-gon, which forces every red edge in the core of every fiber product to also be contained in a red $n$-gon. Therefore the remaining vertices and red edges in the cores each of the fiber products must arise as connected components consisting solely of one red $n$-gon, which is a subgraph of $\sigma(W)$.
    
    If a graph $K\in S$ is $q$-connected, Lemma 3.14 tells us that every connected component of $K\otimes_{U} H$ and $H\otimes_{U} K$ is a subgraph of $K$ for all $H\in S$. While $\beta(K)$ may not already be in $S$ (in which case we add $\beta(K)$ to $S$), $\beta(K)$ is also $q$-connected by Lemma 3.15, and therefore every connected component of $\beta(K)\otimes_{U} H$ and $H\otimes_{U}\beta(K)$ is a subgraph of $\beta(K)$ for all $H\in S$ as well. Since $\beta^2=1$, this guarantees that any $q$-connected $K$ can contribute at most one new graph to $S$, namely $\beta(K)$. In the above calculations we encountered a finite number of $q$-connected relevant graphs, so including these graphs and their images under $\beta$ in $S$ allows $S$ to remain finite for each fixed $n$. Therefore $S$ is finite and satisfies the properties described in Definition 3.7.
\end{proof}

\section{Residual finiteness of $A_{2,3,8}$}

\begin{theorem}
    $A_{2,3,8}$ is residually finite.
\end{theorem}

In this section we will prove Theorem 5.1 in the exact same way that we proved Theorem 4.1. To do so, we will again construct a finite set $S$ that will satisfy Definition 3.7.\\

\begin{proof}
The reason that $A_{2,3,8}$ is a special case comes from the very first calculation, $\sigma(W)\otimes_{U}\sigma(W)$. This calculation is carried out in Figure 7 and results in two connected components. The first connected component is, of course, a copy of $\sigma(W)$ with vertices $(1,1)$, $(2,2)$, $(3,3)$ and $(4,4)$. The second connected component is shown in Figure 28, and will be denoted $Y_1$ and included in $S$.

\begin{figure}[h]
    \centering
    \includegraphics[scale=0.3]{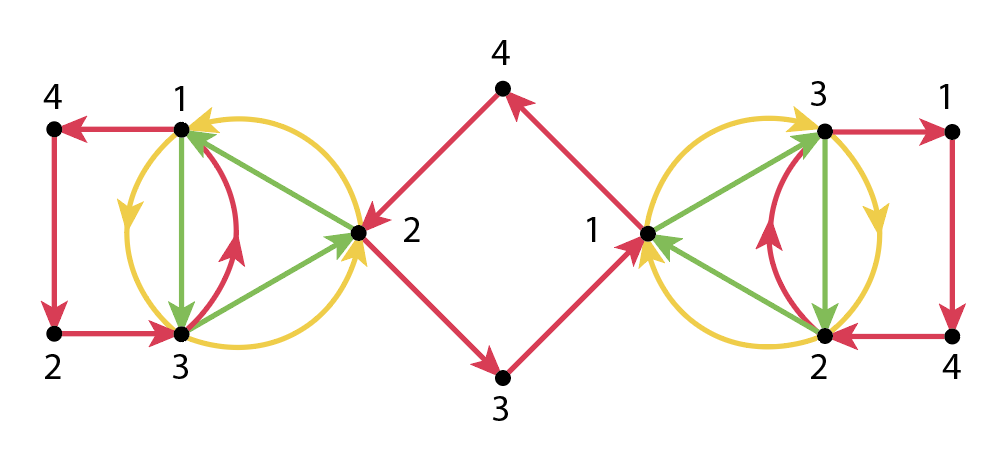}
    \caption{$Y_1$}
    \label{fig:enter-label}
\end{figure}

So far $S=\{\sigma(W),Y_1\}$. The next step is to compute $\beta(Y_1)$. This calculation is performed in Figure 29.

\begin{figure}[h]
    \centering
    \includegraphics[scale=0.35]{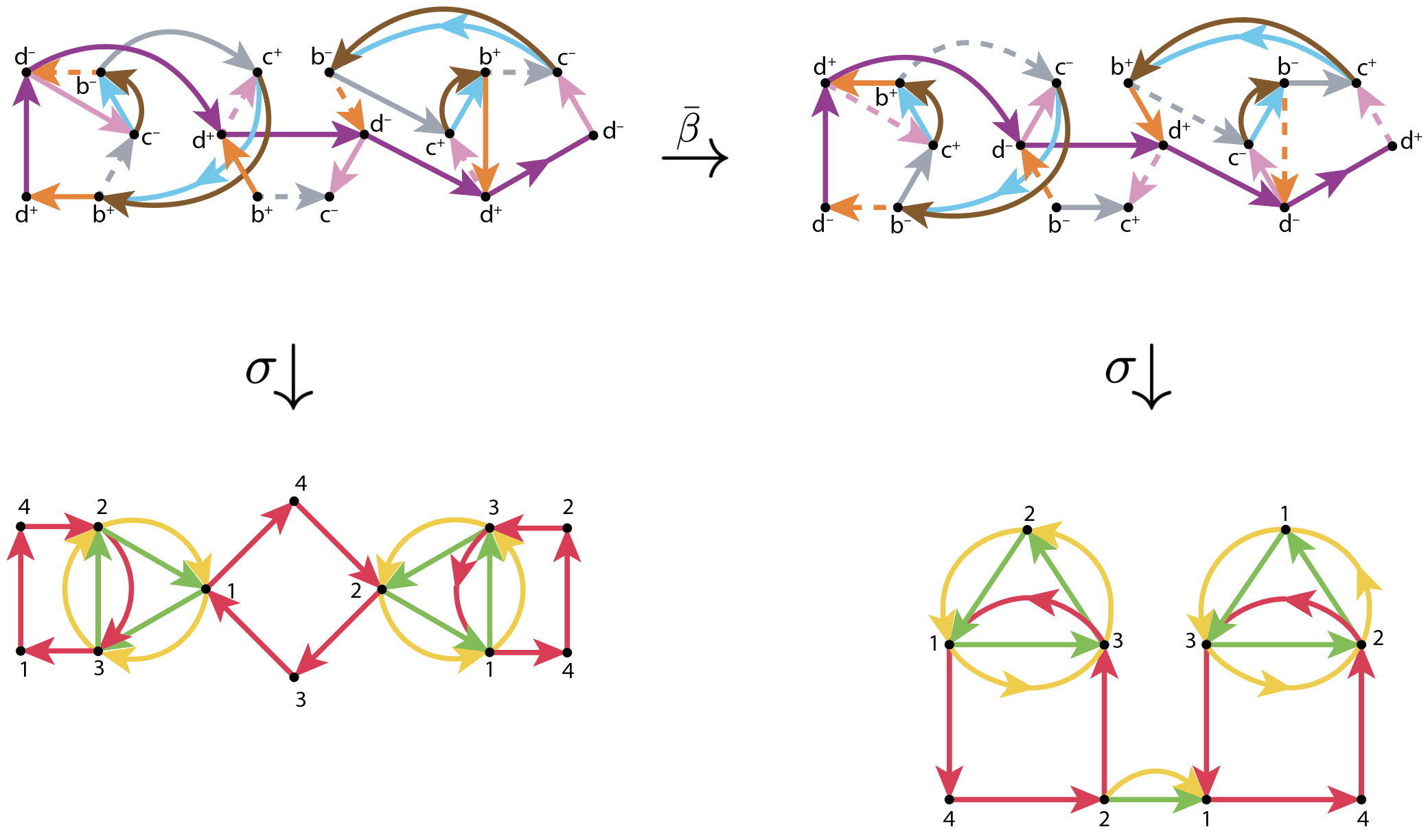}
    \caption{$\beta(Y_1)$}
    \label{fig:enter-label}
\end{figure}

Since $\beta(Y_1)$ is a new graph, we include $\beta(Y_1)$ in $S$ and proceed by performing the fiber product of every pair in $S=\{\sigma(W),Y_1,\beta(Y_1)\}$. These calculations are done in the exact same way as in Section 4, but due to the overwhelming volume of the calculations we have included a Jupyter notebook that can be used to perform the fiber product calculations for us. A full enumeration of all of the connected components that appear in these fiber product computations is shown in Table 1. Readers who would like to perform the fiber product computations should refer to the GitHub link in Section 6. The author used this program to check, by hand, that $\{\sigma(W),Y_1,\beta(Y_1),Y_2,Y_3,\beta(Y_3)\}$ is a complete enumeration of the non $q$-connected components that arise in the computations described by each row in Table 1.

One of the connected components of $\sigma(W)\otimes\beta(Y_1)$ is a new graph that we denote $Y_2$ and is shown in Figure 30. One of the connected components of $\beta(Y_1)\otimes \sigma(W)$ is a new graph that we denote $Y_3$ and is shown in Figure 31.

\begin{figure}[h!]
    \centering
    \includegraphics[scale=0.25]{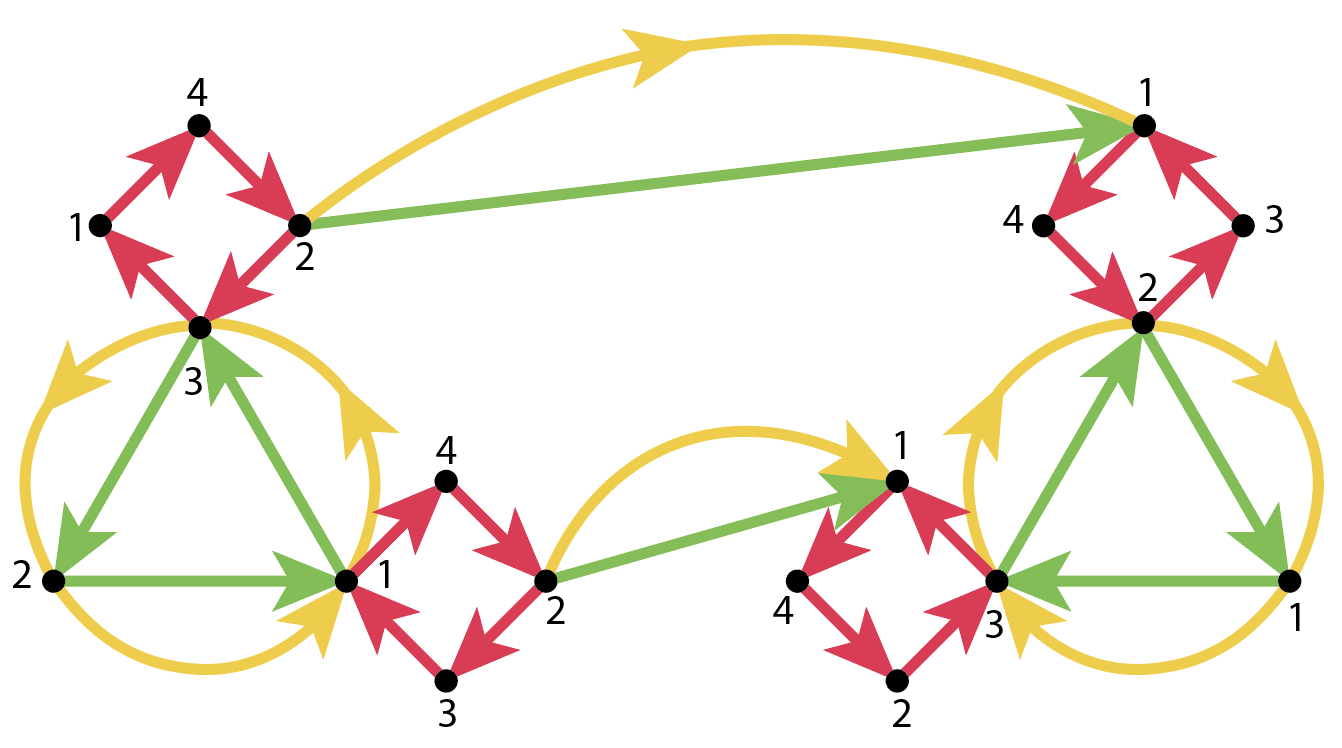}
    \caption{$Y_2$}
    \label{fig:enter-label}
\end{figure}

\begin{figure}[h!]
    \centering
    \includegraphics[scale=0.25]{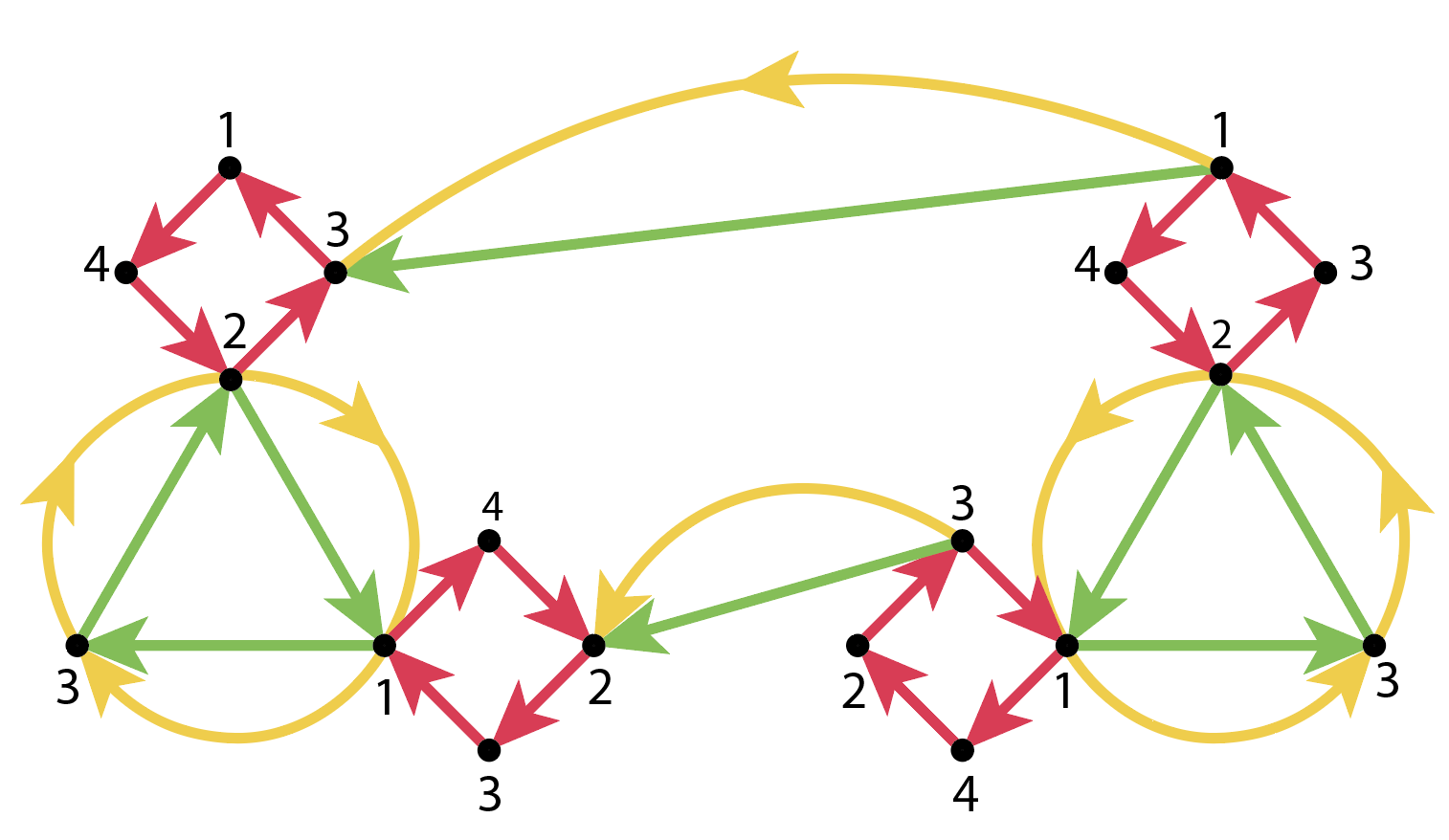}
    \caption{$Y_3$}
    \label{fig:enter-label}
\end{figure}

Again, we need $S$ to be closed under $\beta$, so we must compute $\beta(Y_2)$ and $\beta(Y_3)$. These calculations are shown in Figures 32 \& 33. Luckily, $\beta(Y_2)$ is a rotated copy of $Y_2$, so $\beta(Y_3)$ is the only other new graph that we must add to $S$ after this step. Table 1 contains information about every connected component that arises from the fiber products of every pair of graphs in $S=\{\sigma(W),Y_1,\beta(Y_1),Y_2, Y_3, \beta(Y_3)\}$ where we define $Y_4=Y_2\sqcup Y_3\sqcup\beta(Y_3)$ to streamline computations. Some connected components that arise over the course of these computations are subgraphs of multiple elements in $\{W,Y_1,\beta(Y_1),Y_2,Y_3,\beta(Y_3)\}$. In these situations, a choice was made regarding which column this component is counted in. This choice is arbitrary and does not affect the finiteness of $S$.

\begin{figure}[h!]
    \centering
    \includegraphics[scale=0.3]{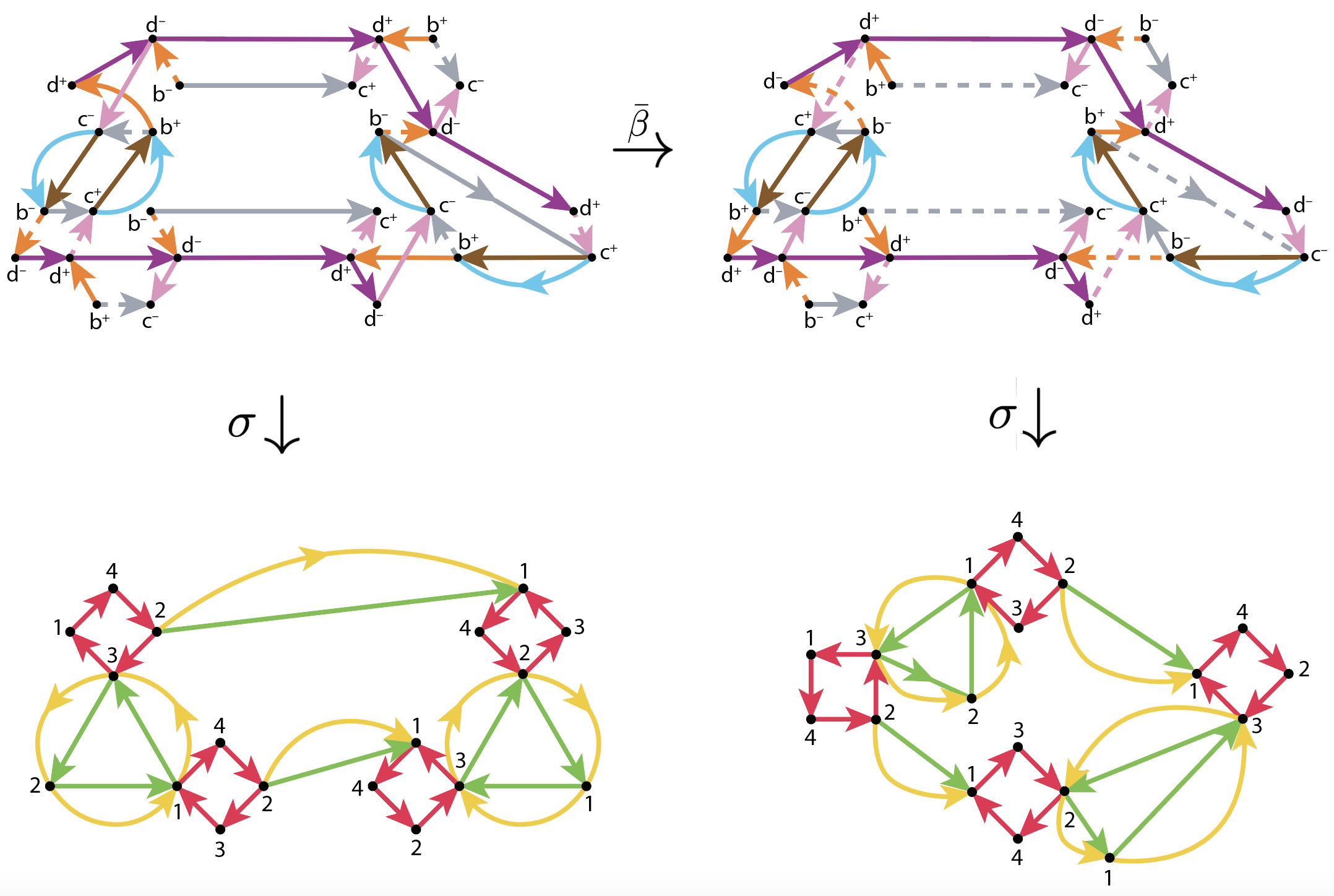}
    \caption{$\beta(Y_2)\cong Y_2$}
    \label{fig:enter-label}
\end{figure}

\begin{figure}[h!]
    \centering
    \includegraphics[scale=0.4]{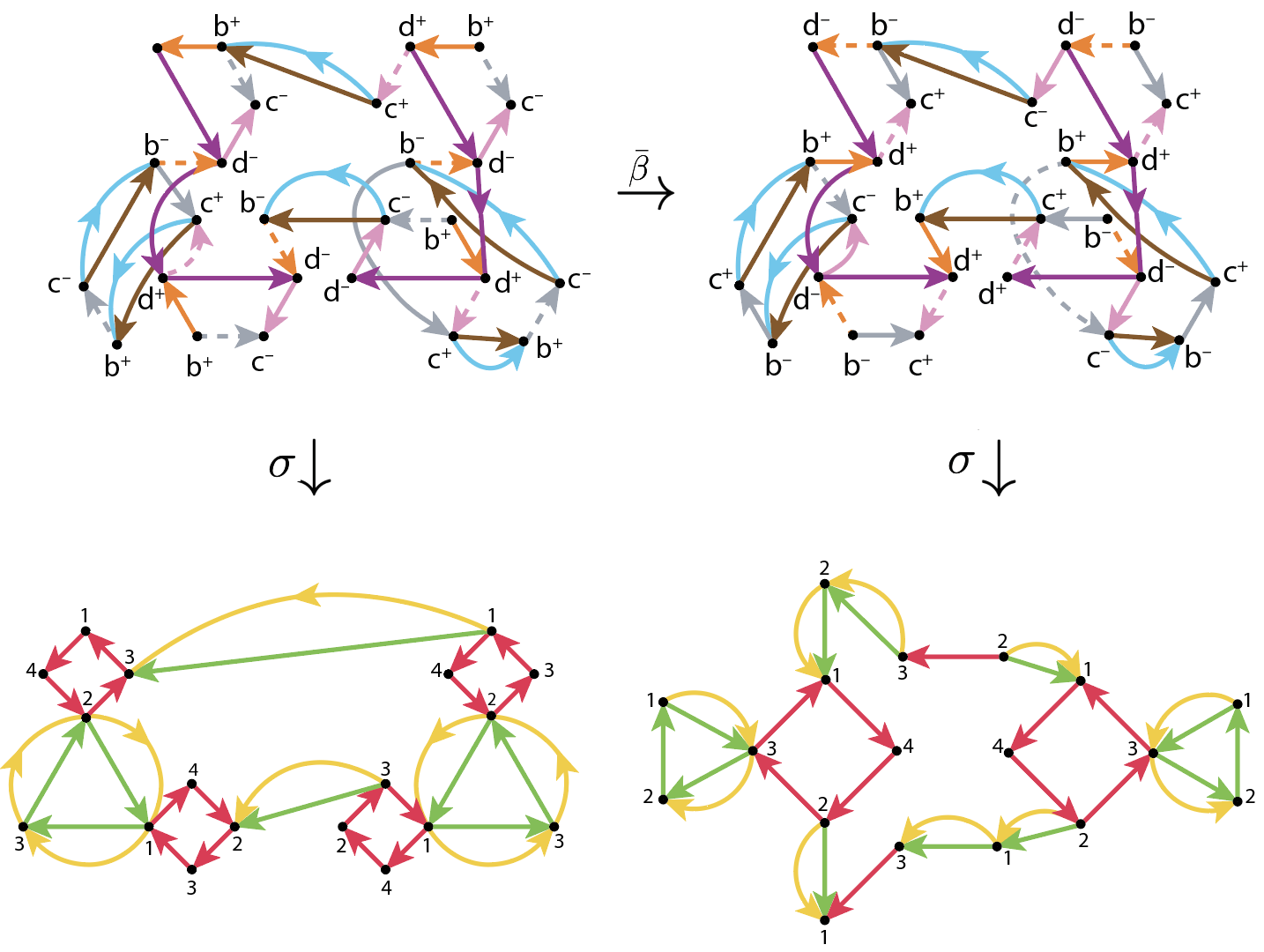}
    \caption{$\beta(Y_3)$}
    \label{fig:enter-label}
\end{figure}

\begin{table}[h!]
\centering
\caption{Number of subgraphs of connected components in $S$}
\label{tab:subgraph-counts}
\begin{tabular}{ |p{2.5cm}||p{1cm}|p{1cm}|p{1cm}|p{1cm}| p{1cm}|p{1cm}|p{2.2cm}|}
 \hline
   & $\sigma(W)$    &$Y_1$&   $\beta(Y_1)$ & $Y_2$ & $Y_3$ & $\beta(Y_3)$ & $q$-connected\\
 \hline
 $\sigma(W)\otimes_{U} \sigma(W)$&1&1&0&0&0&0&0\\
 $\sigma(W)\otimes_{U} Y_1$&0&2&0&0&0&0&2\\
 $\sigma(W)\otimes_{U}\beta(Y_1)$&0&2&1&1&0&0&0\\
 $\sigma(W)\otimes_{U} Y_4$&0&0&0&2&2&2&12\\
 $Y_1\otimes_{U} \sigma(W)$&0&2&0&0&0&0&2\\
 $Y_1\otimes_{U} Y_1$&0&2&0&0&0&0&14\\
 $Y_1\otimes_{U} \beta(Y_1)$&0&4&0&2&0&0&6\\
 $Y_1\otimes_{U} Y_4$&0&0&0&2&2&2&72\\
 $\beta(Y_1)\otimes_{U} \sigma(W)$&0&2&1&0&1&0&0\\
 $\beta(Y_1)\otimes_{U} Y_1$&0&4&0&1&1&0&6\\
 $\beta(Y_1)\otimes_{U}\beta(Y_1)$&0&2&1&1&0&0&8\\
 $\beta(Y_1)\otimes_{U} Y_4$&0&0&0&3&1&1&68\\
 $Y_4\otimes_{U} \sigma(W)$&0&0&0&2&3&1&12\\
 $Y_4\otimes_{U} Y_1$&0&0&0&2&3&1&72\\
 $Y_4\otimes_{U}\beta(Y_1)$&0&0&0&4&0&1&68\\
 $Y_4\otimes_{U} Y_4$&0&0&0&6&3&4&512\\
 \hline
\end{tabular}
\end{table}

Applying the same reasoning as in Section 4.3.5 to Table~\ref{tab:subgraph-counts} proves that $S$ is finite. By Lemma 3.9, this proves that $A_{2,3,8}$ has finite stature with respect to its vertex groups. Therefore, by Theorem 1.3, $A_{2,3,8}$ is residually finite.
\end{proof}

\section{Accessing the Python program}
The program used to analyze the connected components that arise in the $A_{2,3,8}$ fiber products is available at \[\text{https://github.com/GreysonPMeyer/Triangle-Artin-Groups}\]

The program is written in Python, and you
will need a Python interpreter to run it. These interpreters are available, for free and for almost all platforms,
from http://python.org.

\section{Acknowledgements}
The author thanks Kasia Jankiewicz for her guidance \& support, as well as the referee for their insightful feedback.

\bibliographystyle{plain}

\end{document}